\documentclass[final]{siamltex}%
\usepackage{amsmath}
\usepackage{amssymb}
\usepackage{graphicx}
\usepackage{color}
\usepackage{diagrams}
\usepackage{amsfonts}%
\usepackage{stmaryrd}

\newcommand{\enleve}[1]{\relax}
\newcommand{\out}{\mathbin{\otimes}}

\newcommand{\CC}{\mathbb{C}}

\newcommand{\resume}[1]{\relax}

\newcommand{\Parafac}{\textsc{parafac}}
\newcommand{\CanD}{\textsc{CanD}}
\newcommand{\Candecomp}{\textsc{candecomp}}

\newcommand{\bm}[1]{\mathbf{#1}}

\newcommand{\bv}{\bm{v}}

\newcommand{\bx}{\bm{x}}

\newcommand{\bA}{A}
\newcommand{\bB}{B}
\newcommand{\bC}{C}

\newcommand{\bL}{L}
\newcommand{\bM}{M}
\newcommand{\bN}{N}

\newcommand{\cY}{\mathcal{Y}}

\newcommand{\red}[1]{#1}
\newcommand{\blue}[1]{#1}

\newtheorem{example}[theorem]{Example}
\newtheorem{conjecture}[theorem]{Conjecture}

\begin{document}

\title{Symmetric tensors and symmetric tensor rank}
\author{
Pierre Comon\thanks{Laboratoire I3S,
CNRS and the University of Nice, Sophia-Antipolis, France,
\texttt{p.comon@ieee.org}}
\and Gene Golub\thanks{Department of Computer Science and {the} Institute for Computational and Mathematical Engineering, Stanford University, Stanford, CA, USA,
{\texttt{golub}, \texttt{lekheng@cs.stanford.edu}}}
\and Lek-Heng Lim\footnotemark[2]
\and Bernard Mourrain\thanks{Projet GALAAD, INRIA, Sophia-Antipolis, France,
\texttt{mourrain@sophia.inria.fr}}
}
\maketitle

\begin{abstract}
A symmetric tensor is a higher order generalization of a symmetric matrix. In
this paper, we study various properties of symmetric tensors in relation
to a decomposition into a symmetric sum of outer product of vectors.
A rank-$1$ order-$k$ tensor is the outer product of $k$ non-zero vectors.
Any symmetric tensor can be decomposed into a linear combination of rank-$1$ tensors, each of them being
symmetric or not. The \textit{rank} of a symmetric tensor is the minimal number of
rank-$1$ tensors that {is} necessary to reconstruct it.
The \textit{symmetric rank} is obtained when the constituting rank-$1$ tensors are imposed to be themselves
symmetric. It is shown that rank and symmetric rank are equal in a number of
cases, and that they always exist in an algebraically closed field. We will discuss the notion of the generic symmetric rank, which, due to the work of Alexander and Hirschowitz, is now known for any values of dimension and order. We will also show that the set of symmetric tensors of symmetric rank at most $r$ is not closed, unless $r=1$.
\end{abstract}

\begin{keywords} Tensors, multiway arrays, outer product decomposition,
symmetric outer product decomposition, \textsc{candecomp},
\textsc{parafac}, tensor rank, symmetric rank, symmetric tensor
rank, generic {symmetric} rank, {maximal symmetric rank, quantics} \end{keywords}

\begin{AMS} 15A03, 15A21, 15A72, 15A69, 15A18 \end{AMS}

\section{Introduction\label{sec-intro}}

We will be interested in the decomposition of a symmetric tensor into a minimal linear combination of symmetric outer products of vectors (i.e.\ of the form $\mathbf{v}\otimes \mathbf{v}\otimes \dots \otimes \mathbf{v}$). We will see that a decomposition of the form
\begin{equation}\label{eq:sopd}
A = \sum\nolimits_{i=1}^r \lambda_i \mathbf{v}_i\otimes \mathbf{v}_i\otimes \dots \otimes \mathbf{v}_i
\end{equation}
always {exists} for any symmetric tensor $A$ (over any field). One may regard this as a generalization of the eigenvalue decomposition for symmetric matrices to higher order symmetric tensors. In particular, this will allow us to define a notion of \textit{symmetric tensor rank} (as the minimal $r$ over all such decompositions) that {reduces} to the matrix rank for order-$2$ symmetric tensors.

We will call \eqref{eq:sopd} the \textit{symmetric outer product decomposition} of the symmetric tensor $A$ and we will establish its existence in Proposition \ref{symCanD-prop}. This is often abbreviated as \CanD\ in signal processing. The decomposition of a tensor into an (asymmetric) outer product of vectors and the corresponding notion of tensor rank was first introduced and studied by Frank~L.~Hitchcock in 1927 \cite{Hi1, Hi2}. This same decomposition was rediscovered in the 1970s by psychometricians in their attempts to define data analytic models that generalize factor analysis to multiway data \cite{Tuck66:psy}. The name \Candecomp, for `canonical decomposition', was used by Carrol and Chang \cite{CarrC70:psy} while the name \Parafac, for `parallel factor analysis', was used by Harshman \cite{Hars70:ucla} for their respective models.


The symmetric outer product decomposition is particularly important in the process of \textit{blind identification of
under-determined mixtures} (UDM), i.e.\ linear mixtures with more inputs than observable outputs. We refer the reader to \cite{Como02:oxford, ComoR06:SP, DelaDV00:simax2, SidiBG00:ieeesp, SmilBG04} and references therein for a list of other application areas, including speech, mobile communications, machine learning, factor analysis {of $k$-way arrays,} biomedical engineering, psychometrics, and chemometrics.

Despite a growing interest in the symmetric decomposition of symmetric tensors, this topic has not been adequately addressed in the general literature, and even less so in the engineering literature. For several years,
the alternating least squares algorithm has been used to fit data arrays to a
multilinear model \cite{Krus77:laa, SmilBG04}. Yet, the minimization of this
matching error is an ill-posed problem in general, since the set of symmetric tensors of
symmetric rank not more than $r$ is not closed, unless $r=1$ (see Sections \ref{sec-main}
and \ref{sec-examples}) --- a fact that parallels the illposedness discussed in  \cite{dSL}.
The focus of this \red{paper} is mainly on symmetric tensors. The asymmetric case will be addressed in a companion paper, and will use similar tools borrowed from algebraic geometry.

Symmetric tensors {form} a singularly important
class of tensors. Examples where these arise include higher order derivatives
of {smooth functions \cite{L2}}, and {moments and cumulants} of random vectors
\cite{Mccu87}. The decomposition of such symmetric tensors into simpler ones,
as in the symmetric outer product decomposition, plays an
important role in \textit{independent component analysis} \cite{Como02:oxford} and
constitutes a problem of interest {in its own right}. On the other hand the
asymmetric version of the outer product decomposition defined in \eqref{cand:eq} is
central \red{to} \textit{multiway factor analysis} \cite{SmilBG04}.

In Sections \ref{sec-arrays} and \ref{sec-symmetry}, we discuss some classical results in multilinear algebra \cite{B, Gr, L, M, N, Y} and algebraic geometry \cite{Harr98, Z}. While these background materials are well-known to \red{many} pure mathematicians, we found that practitioners and applied mathematicians (in signal processing, neuroimaging, numerical analysis, optimization, etc) --- for whom this paper is intended --- are often unaware of these classical results. For instance, some do not realize that the classical definition of a symmetric tensor given in Definition \ref{defSymmTensor} is equivalent to the requirement that the coordinate array representing the tensor be invariant under all {permutations} of indices, as in Definition \ref{defSymmArray}. Many authors have persistently mislabeled the latter a `\textit{supersymmetric} tensor' (cf.\ \cite{Card91:tor, KofiR02:simax, Q1}). In fact, we have found that even the classical definition of a symmetric tensor is not as well-known as it should be. We see this as an indication of the need to inform our target readership. It is our hope that the background materials presented in Sections \ref{sec-arrays} and \ref{sec-symmetry} will serve such a purpose.

Our contributions will only begin in Section \ref{sec-ranks}, where
the {notions} of maximal and generic rank are analyzed.
The concepts of \textit{symmetry} and \textit{genericity} are recalled in Sections
\ref{sec-symmetry} and \ref{sec-ranks}, respectively. The distinction between
symmetric rank and rank is made in Section \ref{sec-ranks}, and it is
\red{shown} in Section \ref{sec-rank=symrank} \red{that} they \red{must} be
equal {in specific cases}. It is also pointed out in Section \ref{sec-main} that the generic rank
always exists in an algebraically closed field, and that it is not maximal
except in the binary case.
More precisely, the sequence of sets of symmetric tensors of \red{symmetric} rank $r$ increases with $r$ (in the sense of inclusion) up to the generic \red{symmetric} rank, and \textit{decreases} thereafter. In addition, the set of symmetric tensors of \red{symmetric} rank at most $r$ and order $d>2$ is closed only for $r=1$ and $r=R_{\mathsf{S}}$, the maximal \red{symmetric} rank.
Values of the generic \red{symmetric} rank and the {uniqueness} of the {symmetric outer product} decomposition are addressed in Section \ref{sec-values}.
\red{In Section \ref{sec-examples}, we give several examples of sequences of symmetric tensors converging to limits having strictly higher symmetric ranks. We also give an explicit example of a symmetric tensor whose values of symmetric rank over $\mathbb{R}$ and over $\mathbb{C}$ are different.}

In this paper, we restrict our attention \red{mostly} to decompositions {over} the complex
field. \red{A} corresponding study over the real field \red{will require techniques rather different from those} introduced here, as \red{we will elaborate} in Section \ref{realField-sec}.



\section{Arrays and tensors\label{sec-arrays}}

A $k$\textit{-way array} of complex numbers will be written in the form
$\bA =\llbracket a_{j_{1}\cdots j_{k}}\rrbracket_{j_{1},\dots,j_{k}=1}^{n_{1},\dots,n_{k}}$,
where $a_{j_{1}\cdots j_{k}}\in\mathbb{C}$ is the $(j_{1},\dots,j_{k}%
)$-entry of the array. This is sometimes also called a $k$-dimensional \textit{hypermatrix}. We denote the set of all such arrays by $\mathbb{C}^{n_{1}\times\dots\times n_{k}}$, which is evidently a complex vector space of
dimension $n_{1}\cdots n_{k}$ with respect to entry-wise addition and scalar
multiplication. When there is no confusion, we will leave out the range of the
indices and simply write $\bA =\llbracket a_{j_{1}\cdots j_{k}}\rrbracket \in\mathbb{C}%
^{n_{1}\times\dots\times n_{k}}$.

Unless noted otherwise, arrays with at least two indices will be denoted in
uppercase; vectors are one-way arrays, and will be denoted in bold
lowercase. For our purpose, only a few notations related to arrays
\cite{Como02:oxford, DelaDV00:simax2} are necessary.

The \textit{outer product} (or \textit{Segre outer product}) of $k$ vectors
$\mathbf{u}\in\mathbb{C}^{n_{1}},\mathbf{v}\in\mathbb{C}^{n_{2}}%
,\dots,\mathbf{z}\in\mathbb{C}^{n_{k}}$ is defined as
\[
\mathbf{u}\otimes\mathbf{v}\otimes\dots\otimes\mathbf{z} := \llbracket u_{j_{1}}v_{j_{2}%
}\cdots z_{j_{k}}\rrbracket_{j_{1},j_{2},\dots,j_{k}=1}^{n_{1},n_{2},\dots,n_{k}}.
\]
More generally, the outer product of two arrays $\bA $ and $\bB $,
respectively of orders $k$ and $\ell$, is an array of order $k+\ell$,
$\bC =\bA \otimes\bB $ with entries
\[
c_{i_{1}\cdots i_{k}j_{1}\cdots j_{\ell}}:=a_{i_{1}\cdots i_{k}}b_{j_{1}\cdots
j_{\ell}}.
\]
For example, the outer product of two vectors, $\mathbf{u}\otimes\mathbf{v}$,
is a matrix. The outer product of three vectors, or of a matrix with a vector,
is a $3$-way array.

How is an array related to a tensor? Recall that a tensor is simply an element
in the tensor product of vector spaces \cite{B, Gr, L, M, N, Y}. One may
easily check that the so-called Segre map
\begin{align*}
\varphi:\mathbb{C}^{n_{1}}\times\dots\times\mathbb{C}^{n_{k}}  &
\rightarrow\mathbb{C}^{n_{1}\times\dots\times n_{k}},\\
(\mathbf{u},\dots,\mathbf{z})  &  \mapsto\mathbf{u}\otimes\dots\otimes
\mathbf{z}%
\end{align*}
is multilinear. By the universal property of the tensor product \cite{B, Gr,
L, M, N, Y}, there exists a linear map $\theta$
\[
\begin{diagram}[labelstyle=\scriptstyle]
&  & \CC^{n_{1}}\otimes\dots\otimes\CC^{n_{k}}\\
& \ruTo & \dTo^{\theta}\\
\CC^{n_{1}}\times\dots\times\CC^{n_{k}} & \rTo^{\varphi} &
\CC^{n_{1}\times\dots\times n_{k}}.
\end{diagram}
\]
Since $\dim(\mathbb{C}^{n_{1}}\otimes\dots\otimes\mathbb{C}^{n_{k}}%
)=\dim(\mathbb{C}^{n_{1}\times\dots\times n_{k}})$, $\theta$ is an isomorphism
of the vector spaces $\mathbb{C}^{n_{1}}\otimes\dots\otimes\mathbb{C}^{n_{k}}$
and $\mathbb{C}^{n_{1}\times\dots\times n_{k}}$. Consider the canonical basis
of $\mathbb{C}^{n_{1}}\otimes\dots\otimes\mathbb{C}^{n_{k}}$,%
\[
\{\mathbf{e}_{j_{1}}^{(1)}\otimes\dots\otimes\mathbf{e}_{j_{k}}^{(k)}\mid1\leq
j_{1}\leq n_{1},\dots,1\leq j_{k}\leq n_{k}\},
\]
where $\{\mathbf{e}_{1}^{(\ell)},\dots,\mathbf{e}_{n_{\ell}}^{(\ell)}\}$
denotes the canonical basis in $\mathbb{C}^{n_{\ell}}$, $\ell=1,\dots,k$. Then
$\theta$ may be described by%
\[
\theta\Bigl(\sum\nolimits_{j_{1},\dots,j_{k}=1}^{n_{1},\dots,n_{k}}%
a_{j_{1},\dots,j_{k}}\mathbf{e}_{j_{1}}^{(1)}\otimes\dots\otimes
\mathbf{e}_{j_{k}}^{(k)}\Bigr)=\llbracket a_{j_{1}\cdots j_{k}}\rrbracket_{j_{1},\dots,j_{k}%
=1}^{n_{1},\dots,n_{k}}.
\]
So an order-$k$ tensor in $\mathbb{C}^{n_{1}}\otimes\dots\otimes
\mathbb{C}^{n_{k}}$ and a $k$-way array in $\mathbb{C}^{n_{1}\times
\dots\times n_{k}}$ that represents the tensor with respect to a
basis may be regarded as synonymous (up to, of course, the choice of basis). We will
illustrate how the $k$-array representation of an order-$k$ tensor is affected
by a change-of-basis. Let $\bA =\llbracket a_{ijk}\rrbracket \in\mathbb{C}^{n_{1}\times
n_{2}\times n_{3}}$ and let $\bL $, $\bM $, and $\bN $ be
three matrices of size $r_{1}\times n_{1}$, $r_{2}\times n_{2}$, and
$r_{3}\times n_{3}$, respectively. Then \red{the} tensor $\bA $ \red{may be} transformed by
the multilinear map $(\bL ,\bM ,\bN )$ into a tensor
$\bA ^{\prime}= \llbracket a_{pqr}^{\prime}\rrbracket \in\mathbb{C}^{r_{1}\times r_{2}\times
r_{3}}$ defined \red{by}
\begin{equation}
\label{multilin:eq}a_{pqr}^{\prime}=\sum\nolimits_{i,j,k}l_{pi}m_{qj}n_{rk} a_{ijk}.
\end{equation}
When $r_{i}=n_{i}$ and $\bL ,\bM ,\bN $ are nonsingular
matrices\red{, the} above multilinear map may be thought of as a {change-of-bases} (refer to \cite{dSL} for further discussions). \red{We will call this map a \textit{multilinear transform} of $A$.}

In addition to the outer product, we also have an \textit{inner product} or
\textit{contraction product} of two arrays. The mode-$p$ inner product between
two arrays $\bA, \bB $ having the same $p$th dimension is denoted
$\bA {\mathbin{\bullet}}_{p}\bB $, and is obtained by
summing over the $p$th index. More precisely, if $\bA $ and $\bB $
are of orders $k$ and $\ell$ respectively, this yields for $p=1$ the array
$\bC =\bA {\mathbin{\bullet}}_{1}\bB $ of order
$k+\ell-2$:%
\[
c_{i_{_{2}}\cdots i_{k}j_{_{2}}\cdots j_{\ell}}=\sum\nolimits_{\alpha}a_{\alpha
i_{_{2}}\cdots i_{k}}b_{\alpha j_{_{2}}\cdots j_{\ell}}.
\]
Note that some authors \cite{DelaDV00:simax2, ElS, Tuck66:psy} denoted this contraction product as $\bA \times_{p}\bB $ or $\langle \bA, \bB \rangle_p $. By convention, when the contraction is between a tensor and a matrix, it is convenient to assume that the summation is always done on the second matrix index. For instance, the \red{multilinear transform} in \eqref{multilin:eq} may be expressed as $\bA^{\prime}=\bA {\mathbin{\bullet}}_{1}\bL {\mathbin{\bullet}}_{2}\bM {\mathbin{\bullet}}_{3}\bN $. An alternative notation for \eqref{multilin:eq} from the theory of group actions is $A^{\prime} = (L,M,N) \cdot A$, which may be viewed as multiplying $A$ on `three sides' by the matrices $L$, $M$, and $N$ \cite{dSL, J1}.

\section{Symmetric arrays and symmetric tensors}\label{sec-symmetry}

We shall say that a $k$-way array is \textit{cubical} if all its $k$
dimensions are identical, i.e.\ $n_{1}=\dots=n_{k}=n$. A cubical
array will be called \textit{symmetric} if its entries do not change
{under} any permutation of its $k$ indices. Formally, if
$\mathfrak{S}_{k}$ denotes the symmetric group of permutations on
$\{1,\dots,k\}$, then we have

\begin{definition}
\label{defSymmArray}A $k$-way array $\llbracket a_{j_{1}\cdots j_{k}}\rrbracket\in\CC%
^{n\times\dots\times n}$ is called \textbf{symmetric} if
\[
a_{i_{\sigma(1)}\cdots i_{\sigma(k)}}=a_{i_{1}\cdots i_{k}},\qquad i_{1}%
,\dots,i_{k}\in\{1,\dots,n\},
\]
for all permutations $\sigma\in\mathfrak{S}_{k}$.
\end{definition}

For example, a $3$-way array $\llbracket a_{ijk} \rrbracket \in
\mathbb{C}^{n \times n \times n}$ is symmetric if
\[
a_{ijk} = a_{ikj} = a_{jik} = a_{jki} = a_{kij} = a_{kji}
\]
\red{for all $i,j,k \in \{1,\dots,n\}$.}

Such arrays have been improperly labeled `supersymmetric'\ tensors
(cf.\ \cite{Card91:tor, KofiR02:simax, Q1} among others); this
terminology should be avoided since it refers to an entirely different class
of tensors \cite{BrinHT92:amath}. The word `supersymmetric' has
\textit{always} been used in both mathematics and physics \cite{F, V, W} to
describe objects with a $\mathbb{Z}_{2}$-\red{grading and} so using it in the
sense of \cite{Card91:tor, KofiR02:simax, Q1} is both inconsistent and confusing
(the correct usage will be one in the sense of \cite{BrinHT92:amath}). In
fact, we will show below \red{in Proposition \ref{propSame}} that there is no
difference between Definition \ref{defSymmArray} and the usual definition of a
symmetric tensor in mathematics \cite{B, Gr, L, M, N, Y}. In other words, the
prefix `super' in `supersymmetric tensor'\red{, when used} in the sense of \cite{Card91:tor,
KofiR02:simax, Q1}\red{,} is superfluous.

We will write {$\mathsf{T}^{k}(\CC^{n}):=\CC^{n}\otimes
\dots\otimes\CC^{n}$ ($k$ copies)}, the set of all order-$k$ dimension-$n$ cubical tensors.
We define a group action $\mathfrak{S}_{k}$ on $\mathsf{T}^{k}%
(\CC^{n})$ as follows. For any $\sigma\in\mathfrak{S}_{k}$ and
$\bx_{i_{1}}\otimes\dots\otimes\bx_{i_{k}}\in\mathsf{T}^{k}(\CC^{n})$, we let
\[
\sigma(\bx_{i_{1}}\otimes\dots\otimes\bx_{i_{k}}):=
\bx_{i_{\sigma(1)}}\otimes\dots\otimes\bx_{i_{\sigma(k)}}
\]
and extend this linearly to all of $\mathsf{T}^{k}(\CC^{n})$. Thus each
$\sigma\in\mathfrak{S}_{k}$ defines a linear operator $\sigma:\mathsf{T}%
^{k}(\CC^{n})\rightarrow\mathsf{T}^{k}(\CC^{n})$. The standard
definition of a symmetric tensor in mathematics \cite{B, Gr, L, M, N, Y} looks
somewhat different from Definition \ref{defSymmArray} and is given as follows.

\begin{definition}
\label{defSymmTensor}An order-$k$ tensor $\bA\in\mathsf{T}^{k}(\CC^{n})$ is \textbf{symmetric} if
\begin{equation}
\sigma(\bA)=\bA \label{eqLO}
\end{equation}
for all permutations $\sigma\in\mathfrak{S}_{k}$. The set of symmetric tensors
in $\mathsf{T}^{k}(\CC^{n})$ will be denoted by $\mathsf{S}^{k}(\CC^{n})$.
\end{definition}

Let $S:\mathsf{T}^{k}(\CC^{n}) \rightarrow \mathsf{T}^{k}(\CC^{n})$ be the linear operator defined by%
\[
S:=\frac{1}{k!}\sum\nolimits_{\sigma\in\mathfrak{S}_{k}}\sigma.
\]
Note that given any $\sigma\in\mathfrak{S}_{k}$,
\[
\sigma\circ S=S\circ\sigma=S.
\]
Here $\circ$ denotes the composition of the linear operators $\sigma$ and $S$.

\begin{proposition}\label{prop:symm-equiv}
An order-$k$ tensor $\bA\in\mathsf{T}^{k}(\CC^{n})$ is symmetric
if and only if
\[
S(\bA):=\frac{1}{k!}\sum\nolimits_{\sigma\in\mathfrak{S}_{k}}\sigma
(\bA)=\bA.
\]

\end{proposition}

\begin{proof}
Clearly, if $\bA$ is symmetric, then%
\[
S(\bA)=\frac{1}{k!}\sum\nolimits_{\sigma\in\mathfrak{S}_{k}}\sigma
(\bA)=\frac{1}{k!}\sum\nolimits_{\sigma\in\mathfrak{S}_{k}}\bA%
=\bA.
\]
Conversely, if $S(\bA)=\bA$, then%
\[
\sigma(\bA)=\sigma(S(\bA))=\sigma\circ S(\bA%
)=S(\bA)=\bA%
\]
for all $\sigma\in\mathfrak{S}_{k}$; and so $\bA$ is symmetric.
\end{proof}

In other words, a symmetric tensor is an eigenvector of the linear
operator $S$ {with} eigenvalue $1$.
$\mathsf{S}^{k}(\CC^{n})$ is the $1$-eigenspace of
$S:\mathsf{T}^{k}(\CC^{n})\rightarrow \mathsf{T}^{k}(\CC^{n})$. {Proposition \ref{prop:symm-equiv} implies} that
$\mathsf{S}^{k}(\CC^{n})=S(\mathsf{T}^{k}(\CC^{n}))$ {and it} is also
easy to see that $S$ is a projection of $\mathsf{T}^{k}(\CC^{n})$
onto {the subspace} $\mathsf{S}^{k}(\CC^{n})$, i.e.\ $S^{2}=S$.

\subsection{Equivalence with homogeneous polynomials\label{poly:subsec}}

We adopt the following standard {shorthand. For} any $\mathbf{e}_{i_1}, \dots, \mathbf{e}_{i_k}\in\mathbb{C}^{n}$ { with $i_1,\dots, i_k \in \{ 1,\dots, n\}$, we write}
\begin{equation}
\mathbf{e}_{i_{1}}\dotsb\mathbf{e}_{i_{k}}:=S(\mathbf{e}_{i_{1}%
}\otimes\dotsb\otimes\mathbf{e}_{i_{k}})=\frac{1}{k!}\sum\nolimits_{\sigma
\in\mathfrak{S}_{k}}\mathbf{e}_{i_{\sigma(1)}}\otimes\dotsb\otimes
\mathbf{e}_{i_{\sigma(k)}}. \label{sh1}%
\end{equation}
Then since $S\sigma=S$, the term $\mathbf{e}_{i_{1}}%
\dotsb\mathbf{e}_{i_{k}}$ depends only on the number of times each
$\mathbf{e}_{i}$ enters this product and we may write
\begin{equation}
\mathbf{e}_{i_{1}}\dotsb\mathbf{e}_{i_{k}}=\mathbf{e}_{1}^{p_{1}}%
\dotsb\mathbf{e}_{n}^{p_{n}} \label{sh2}%
\end{equation}
where $p_{i}$ is the multiplicity (which may be $0$) of occurrence of
$\mathbf{e}_{i}$ in $\mathbf{e}_{i_{1}}\dotsb\mathbf{e}_{i_{k}}$. Note that
$p_{1},\dots,p_{n}$ are nonnegative integers satisfying $p_{1}+\dotsb
+p_{n}=k$.

\begin{proposition}
\label{propBasis}Let $\{\mathbf{e}_{1},\dotsc,\mathbf{e}_{n}\}$ be a basis of
$\mathbb{C}^{n}$. Then
\[
\{S(\mathbf{e}_{i_{1}}\otimes\dotsb\otimes\mathbf{e}_{i_{k}})\mid1\leq
i_{1}\leq\dotsb\leq i_{k}\leq n\}
\]
or, explicitly,%
\[
\Bigl\{
\frac{1}{k!}\sum\nolimits_{\sigma\in\mathfrak{S}_{k}}\mathbf{e}_{i_{\sigma
(1)}}\otimes\dotsb\otimes\mathbf{e}_{i_{\sigma(k)}}\Bigm| 1\leq i_{1}%
\leq\dotsb\leq i_{k}\leq n \Bigr\},
\]
is a basis of $\mathsf{S}^{k}(\mathbb{C}^{n})$. Furthermore,
\[
\dim_{\mathbb{C}}\mathsf{S}^{k}(\mathbb{C}^{n})=\dbinom{n+k-1}{k}.
\]

\end{proposition}

\begin{proof}
Since \red{$\mathcal{B} = \{\mathbf{e}_{i_{1}}\otimes\dotsb\otimes\mathbf{e}_{i_{k}}\mid1\leq
i_{1}\leq n,\dots,1\leq i_{k}\leq n\}$} is a basis for $\mathsf{T}%
^{k}(\mathbb{C}^{n})$ and since $S$ maps $\mathsf{T}^{k}(\mathbb{C}^{n})$ onto
$\mathsf{S}^{k}(\mathbb{C}^{n})$, the set%
\[
S(\mathcal{B})=\{\mathbf{e}_{i_{1}}\dotsb\mathbf{e}_{i_{k}}\mid1\leq i_{1}%
\leq\dotsb\leq i_{k}\leq n\}\newline=\{\mathbf{e}_{1}^{p_{1}}\dotsb
\mathbf{e}_{n}^{p_{n}}\mid p_{1}+\dotsb+p_{n}=k\}
\]
spans $\mathsf{S}^{k}(\mathbb{C}^{n})$. Vectors in $S(\mathcal{B})$ are
linearly independent: if $(p_{1},\dotsc,p_{n})\neq(q_{1},\dotsc,q_{n})$, then
the tensors $\mathbf{e}_{1}^{p_{1}}\dotsb\mathbf{e}_{n}^{p_{n}}$ and
$\mathbf{e}_{1}^{q_{1}}\dotsb\mathbf{e}_{n}^{q_{n}}$ are respectively linear
combinations of two non-intersecting subsets of basis elements of
$\mathsf{T}^{k}(\mathbb{C}^{n})$. The cardinality of $S(\mathcal{B})$ is
precisely number of partitions of $k$ into a sum of $n$ nonnegative integers,
i.e.\ $\binom{n+k-1}{k}$.
\end{proof}

If we regard $\mathbf{e}_{j}$ in \eqref{sh2} as variables (i.e.\ indeterminates), then every
symmetric tensor of order $k$ and dimension $n$ may be uniquely associated
with a \textit{homogeneous polynomial} of degree $k$ in $n$ variables. Recall that these are just polynomials in $n$ variables whose constituting monomials all have the same total degree $k$.
Homogeneous polynomials are also called \textit{quantics} and those of degrees $1$, $2$, and $3$ are often called \textit{linear forms}, \textit{quadratic forms}, and \textit{cubic forms} (or just \textit{cubics}) respectively.
From now on, we will use more standard notation for the variables --- $x_j$ instead of $\mathbf{e}_{j}$. So the monomial on the \textsc{rhs} of \eqref{sh2} now becomes $x_{1}^{p_{1}}\dotsb x_{n}^{p_{n}}$.
To further simplify this notation, we will adopt the following standard \textit{multi-index notations}:
\[
\mathbf{x}^{\boldsymbol{p}}:=\prod\nolimits_{k=1}^{n}x_{k}^{p_{k}} \quad \text{and} \quad \lvert
\boldsymbol{p}\rvert:=\sum\nolimits_{k=1}^{n}p_{k},
\]
where {$\boldsymbol{p}$ denotes a $k$-vector of nonnegative integers}. We will also write $\mathbb{C}[x_{1},\dots,x_{n}]_{k}$ for the set of homogeneous polynomials of degree $k$ in $n$ variables (again a
standard notation). Then any symmetric tensor $\llbracket a_{j_{1}\cdots j_{k}%
}\rrbracket = \llbracket a_{\boldsymbol{j}}\rrbracket \in\mathsf{S}^{k}(\mathbb{C}^{n})$ can be associated with a
unique homogeneous polynomial $F\in$ $\mathbb{C}[x_{1},\dots,x_{n}]_{k}$ via
the \red{expression}
\begin{equation}
F(\mathbf{x})=\sum\nolimits_{\boldsymbol{j}}a_{\boldsymbol{j}}\mathbf{x}^{\boldsymbol{p}%
(\boldsymbol{j})}, \label{polytens:eq}%
\end{equation}
where for every $\boldsymbol{j}=(j_1,\dots,j_k)$, one associates
bijectively the nonnegative integer vector $\boldsymbol{p}(\boldsymbol{j}) = (p_1(\boldsymbol{j}),\dots,p_n(\boldsymbol{j}))$ with
$p_{j}(\boldsymbol{j})$ counting the number of times
index $j$ appears in $\boldsymbol{j}$ \red{\cite{ComoM96:SP, Como02:oxford}}. We have in particular $\lvert
\boldsymbol{p}(\boldsymbol{j})\rvert=k$. The converse is true as well, and the
correspondence between symmetric tensors and homogeneous polynomials is
obviously bijective. Thus
\begin{equation}\label{eq:symm-hpoly}
\mathsf{S}^{k}(\mathbb{C}^{n})\cong\mathbb{C}[x_{1},\dots,x_{n}]_{k}.
\end{equation}

This justifies the use of the Zariski topology, where {the elementary closed
subsets are the} common zeros of a finite number of homogeneous polynomials
\cite{Shaf77}. Note that for asymmetric tensors, the same association is not
possible (although they can still be associated with polynomials via another
bijection). As will be subsequently seen, this identification of symmetric
tensors with homogeneous polynomials will allow us to prove some interesting
facts about symmetric tensor rank.

We {will now proceed to} define a {useful} `inner product' on $\mathbb{C}[x_{1},\dots,x_{n}]_{k}$. {For any} $F,G\in\mathbb{C}[x_{1},\dots,x_{n}]_{k}$ {written} as%
\[
F(\mathbf{x})=\sum\nolimits_{\lvert\boldsymbol{p}\rvert=k}\binom{k}{p_{1},\dots,p_{n}%
}a_{\boldsymbol{p}}\mathbf{x}^{\boldsymbol{p}},\qquad G(\mathbf{x})=\sum\nolimits
_{\lvert\boldsymbol{p}\rvert=k}\binom{k}{p_{1},\dots,p_{n}}b_{\boldsymbol{p}%
}\mathbf{x}^{\boldsymbol{p}},
\]
we let%
\[
\langle F,G\rangle:=\sum\nolimits_{\lvert\boldsymbol{p}\rvert=k}\binom{k}{p_{1},\dots
,p_{n}}a_{\boldsymbol{p}}b_{\boldsymbol{p}}=\sum\nolimits_{p_{1}+\dots+p_{n}=k}\frac{k!}%
{p_{1}!\cdots p_{n}!}a_{p_{1}\cdots p_{n}}b_{p_{1}\cdots p_{n}}.
\]
Note that $\langle\cdot,\cdot\rangle$ cannot be an inner product in the usual
sense since $\langle F,F\rangle$ is in general complex valued (recall that for
an inner product, we will need $\langle F,F\rangle\geq0$ for all $F$).
However, we will show that it is a non-degenerate symmetric bilinear form.

\begin{lemma}
\label{lemNondegen}The bilinear form $\langle\cdot,\cdot\rangle:\mathbb{C}%
[x_{1},\dots,x_{n}]_{k}\times\mathbb{C}[x_{1},\dots,x_{n}]_{k}\rightarrow
\mathbb{C}$ defined above is symmetric and non-degenerate. \red{In other words,} $\langle
F,G\rangle=\langle G,F\rangle$ for \red{every} $F,G\in\mathbb{C}[x_{1},\dots
,x_{n}]_{k}$; and if $\langle F,G\rangle=0$ for all $G\in\mathbb{C}%
[x_{1},\dots,x_{n}]_{k}$, then $F\equiv0$.
\end{lemma}

\begin{proof}
The bilinearity and symmetry is immediate from definition. Suppose $\langle
F,G\rangle=0$ for all $G\in\mathbb{C}[x_{1},\dots,x_{n}]_{k}$. Choose $G$ to
be the monomials%
\[
G_{\boldsymbol{p}}(\mathbf{x})=\binom{k}{p_{1},\dots,p_{n}}\mathbf{x}^{\boldsymbol{p}}%
\]
where $\lvert\boldsymbol{p}\rvert=k$ and we see immediately that%
\[
a_{\boldsymbol{p}}=\langle F,G_{\boldsymbol{p}}\rangle=0.
\]
Thus $F\equiv0$.
\end{proof}

In the special case where $G$ is the $k$th power of a linear form, we have the following lemma.
The main interest in introducing this inner product lies precisely in establishing this lemma.

\begin{lemma}
\label{lemDual}Let $G=(\beta_{1}x_{1}+\dots+\beta_{n}x_{n})^{k}$. Then for any
$F\in\mathbb{C}[x_{1},\dots,x_{n}]_{k}$, we have%
\[
\langle F,G\rangle=F(\beta_{1},\dots,\beta_{n}),
\]
i.e.\ $F$ evaluated at {$(\beta_{1},\dots,\beta_{n}) \in\mathbb{C}^{n}$}.
\end{lemma}

\begin{proof}
Let $b_{\boldsymbol{p}}=\beta_{1}^{p_{1}}\cdots\beta_{n}^{p_{n}}$ for all
$\boldsymbol{p}=(p_{1},\dots,p_{n})$ such that $\lvert\boldsymbol{p}\rvert=k$. The
multinomial expansion then yields%
\[
(\beta_{1}x_{1}+\dots+\beta_{n}x_{n})^{k}=\sum\nolimits_{\lvert\boldsymbol{p}\rvert
=k}\binom{k}{p_{1},\dots,p_{n}}b_{\boldsymbol{p}}\mathbf{x}^{\boldsymbol{p}}.
\]
For any $F(\mathbf{x})=\sum\nolimits_{\lvert\boldsymbol{p}\rvert=k}\binom{k}{p_{1}%
,\dots,p_{n}}a_{\boldsymbol{p}}\mathbf{x}^{\boldsymbol{p}}$,
\[
F(\beta_{1},\dots,\beta_{n})=\sum\nolimits_{\lvert\boldsymbol{p}\rvert=k}\binom{k}%
{p_{1},\dots,p_{n}}a_{\boldsymbol{p}}b_{\boldsymbol{p}}=\langle F,G\rangle
\]
as required.
\end{proof}

\subsection{Equivalence with usual definition}

As mentioned earlier, we will show that a tensor is symmetric in the sense of
Definition \ref{defSymmTensor} if and only if its corresponding array is
symmetric in the sense of Definition \ref{defSymmArray}.

\begin{proposition}
\label{propSame}Let $\bA \in\mathsf{T}^{k}(\mathbb{C}^{n})$ and
$\llbracket a_{j_{1}\cdots j_{k}}\rrbracket\in\mathbb{C}^{n\times\dots\times n}$ be its
corresponding $k$-array. Then%
\[
\sigma(\bA )=\bA %
\]
for all permutations $\sigma\in\mathfrak{S}_{k}$ if and only if%
\[
a_{i_{\sigma(1)}\cdots i_{\sigma(k)}}=a_{i_{1}\cdots i_{k}},\qquad i_{1}%
,\dots,i_{k}\in\{1,\dots,n\},
\]
for all permutations $\sigma\in\mathfrak{S}_{k}$.
\end{proposition}

\begin{proof}
Suppose $\llbracket a_{i_{1}\cdots i_{k}}\rrbracket \in\mathbb{C}^{n\times\dots\times n}$ is
symmetric in the sense of Definition \ref{defSymmArray}. Then the
corresponding tensor%
\[
\bA =\sum\nolimits_{i_{1},\dots,i_{k}=1}^{n}a_{i_{1}\cdots i_{k}}\mathbf{e}%
_{i_{1}}\otimes\dots\otimes\mathbf{e}_{i_{k}},
\]
where $\{\mathbf{e}_{1},\dots,\mathbf{e}_{n}\}$ denotes the canonical basis in
$\mathbb{C}^{n}$, satisfies the following:%
\begin{align*}
S(\bA )  &  =\sum\nolimits_{i_{1},\dots,i_{k}=1}^{n}a_{i_{1}\cdots i_{k}%
}S(\mathbf{e}_{i_{1}}\otimes\dots\otimes\mathbf{e}_{i_{k}}) &  &  \text{($S$
linear)}\\
&  =\frac{1}{k!}\sum\nolimits_{i_{1},\dots,i_{k}=1}^{n}a_{i_{1}\cdots i_{k}}\left[
\sum\nolimits_{\sigma\in\mathfrak{S}_{k}}\mathbf{e}_{i_{\sigma(1)}}\otimes\dots
\otimes\mathbf{e}_{i_{\sigma(k)}}\right]  &  & \\
&  =\frac{1}{k!}\sum\nolimits_{i_{1},\dots,i_{k}=1}^{n}\left[  \sum\nolimits_{\sigma
\in\mathfrak{S}_{k}}a_{i_{\sigma(1)}\cdots i_{\sigma(k)}}\right]
\mathbf{e}_{i_{1}}\otimes\dots\otimes\mathbf{e}_{i_{k}} &  & \\
&  =\frac{1}{k!}\sum\nolimits_{i_{1},\dots,i_{k}=1}^{n}k!a_{i_{1}\cdots i_{k}%
}\mathbf{e}_{i_{1}}\otimes\dots\otimes\mathbf{e}_{i_{k}} &  &
\text{($\llbracket a_{i_{1}\cdots i_{k}}\rrbracket$ symmetric)}\\
&  =\bA . &  &
\end{align*}
Hence $\bA $ is a symmetric tensor in the sense of Definition
\ref{defSymmTensor}.

Conversely, let $\bA $ $\in\mathsf{T}^{k}(\mathbb{C}^{n})$ be symmetric
in the sense of Definition \ref{defSymmTensor} and%
\[
\bA =\sum\nolimits_{i_{1},\dots,i_{k}=1}^{n}a_{i_{1}\cdots i_{k}}\mathbf{e}%
_{i_{1}}\otimes\dots\otimes\mathbf{e}_{i_{k}}%
\]
be the expression of $\bA $ with respect to $\{\mathbf{e}_{i_{1}}%
\otimes\dots\otimes\mathbf{e}_{i_{k}}\mid1\leq i_{1},\dots i_{k}\leq n\}$, the
canonical basis of $\mathsf{T}^{k}(\mathbb{C}^{n})$. Then%
\[
S(\bA )=\bA %
\]
implies%
\[
\sum\nolimits_{i_{1},\dots,i_{k}=1}^{n}\left[  \frac{1}{k!}\sum\nolimits_{\sigma\in
\mathfrak{S}_{k}}a_{i_{\sigma(1)}\cdots i_{\sigma(k)}}\right]  \mathbf{e}%
_{i_{1}}\otimes\dots\otimes\mathbf{e}_{i_{k}}=\sum\nolimits_{i_{1},\dots,i_{k}=1}%
^{n}a_{i_{1}\cdots i_{k}}\mathbf{e}_{i_{1}}\otimes\dots\otimes\mathbf{e}%
_{i_{k}}.
\]
Since $\{\mathbf{e}_{i_{1}}\otimes\dots\otimes\mathbf{e}_{i_{k}}\mid1\leq
i_{1},\dots i_{k}\leq n\}$ is a linearly independent set, we must have%
\begin{equation}
\frac{1}{k!}\sum\nolimits_{\sigma\in\mathfrak{S}_{k}}a_{i_{\sigma(1)}\cdots
i_{\sigma(k)}}=a_{i_{1}\cdots i_{k}}\qquad\text{for all }i_{1},\dots,i_{k}%
\in\{1,\dots,n\}. \label{av}%
\end{equation}
For any given $\tau\in\mathfrak{S}_{k}$, we have
\begin{align*}
a_{i_{\tau(1)}\cdots i_{\tau(k)}}  &  =\frac{1}{k!}\sum\nolimits_{\sigma\in
\mathfrak{S}_{k}}a_{i_{\sigma(\tau(1))}\cdots i_{\sigma(\tau(k))}} &  &
\text{(by \eqref{av})}\\
&  =\frac{1}{k!}\sum\nolimits_{\sigma\in\tau\mathfrak{S}_{k}}a_{i_{\sigma(1)}\cdots
i_{\sigma(k)}} &  & \\
&  =\frac{1}{k!}\sum\nolimits_{\sigma\in\mathfrak{S}_{k}}a_{i_{\sigma(1)}\cdots
i_{\sigma(k)}} &  &  \text{($\tau\mathfrak{S}_{k}=\mathfrak{S}_{k}$ as
$\mathfrak{S}_{k}$ is a group)}\\
&  =a_{i_{1}\cdots i_{k}} &  &  \text{(by \eqref{av}).}%
\end{align*}
Since this holds for arbitrary $\tau\in\mathfrak{S}_{k}$, the array
$\llbracket a_{i_{1}\cdots i_{k}}\rrbracket$ is symmetric in the sense of Definition
\ref{defSymmArray}.
\end{proof}


\section{Notions of rank for symmetric tensors}\label{sec-ranks}

We will discuss two notions of rank for symmetric tensors --- the outer product rank (defined for all tensors) and the symmetric outer product rank (defined only for symmetric tensors). We will show that under certain conditions, they are one and the same. However it is not known if they are equal on all symmetric tensors \red{in general}.

\subsection{Outer product decomposition and rank}

Any tensor can always be decomposed (possibly non-uniquely) as:
\begin{equation}
\bA =\sum\nolimits_{i=1}^{r}\mathbf{u}_{i}\otimes\mathbf{v}_{i}\otimes
\cdots\otimes\mathbf{w}_{i}. \label{cand:eq}%
\end{equation}
The \textit{tensor rank}, $\operatorname*{rank}(\bA )$, is defined as
the smallest integer $r$ such that this decomposition holds exactly \cite{Hi1,
Hi2}. Among other properties, note that this outer product decomposition
remains valid in a ring, and that an outer product decomposition of a
multilinear transform of $A$ equals the multilinear transform of an
outer product decomposition of $\bA $. In other words, if \eqref{cand:eq} is an outer product decomposition of $\bA$, then
\[
\bA {\mathbin{\bullet}} _{1}\bL {\mathbin{\bullet}}_{2}\bM {\mathbin{\bullet}}_{3}\cdots{\mathbin{\bullet}}_{k}\bN =\sum\nolimits_{i=1}^{r}\bL \mathbf{u}_{i} \otimes \bM \mathbf{v}_{i} \otimes\cdots\otimes \bN\mathbf{w}_{i}
\]
is an outer product decomposition of {$\bA {\mathbin{\bullet}} _{1}\bL {\mathbin{\bullet}}_{2}\bM {\mathbin{\bullet}}_{3}\cdots{\mathbin{\bullet}}_{k}\bN $, which may also be written as $(L,M,\dots, N)\cdot A$}. The outer product decomposition has often been regarded synonymously as the data analytic models \textsc{candecomp} \cite{CarrC70:psy} and \textsc{parafac} \cite{Hars70:ucla} {where} the decomposition is used to analyze multiway psychometric data.

\begin{definition}
The \textit{rank} of $\bA =\llbracket a_{j_{1}\cdots j_{k}}\rrbracket \in\mathbb{C}%
^{d_{1}\times\dots\times d_{k}}$ is defined as%
\[
\operatorname*{rank}(\bA ):=\min\{r\mid\bA ={\textstyle\sum
\nolimits_{i=1}^{r}} \mathbf{u}_{i}\otimes\mathbf{v}_{i}\otimes\dots
\otimes \mathbf{w}_{i}\}.
\]
If $\bA =\llbracket a_{j_{1}\cdots j_{k}}\rrbracket \in\mathsf{S}^{k}(\mathbb{C}^{n})$, then
we may also define the notion of \textit{symmetric rank} via
\[
\operatorname*{rank}\nolimits_{\mathsf{S}}(\bA ):=\min\{s\mid
\bA ={\textstyle\sum\nolimits_{i=1}^{s}} \mathbf{y}_{i}\otimes
\dots\otimes\mathbf{y}_{i}\}.
\]

\end{definition}

Note that over $\mathbb{C}$, the coefficients $\lambda_i$ appearing in decomposition \eqref{eq:sopd} may be set to $1$; this is legitimate since any {complex} number admits a $k$th root in $\mathbb{C}$. Henceforth, we will adopt the following notation
\begin{equation}
\mathbf{y}^{\otimes k}:=\overbrace{\mathbf{y}
\otimes\dots\otimes\mathbf{y}}^{k\text{ copies}}.\label{sh3}%
\end{equation}

If in \eqref{cand:eq}, we have $\mathbf{u}_{i}=\mathbf{v}_{i}=\dots
=\mathbf{w}_{i}$ for every $i$, then we may call it a \red{\textit{symmetric
outer product decomposition}}, yielding a \textit{symmetric rank}, $\operatorname*{rank}
\nolimits_{\mathsf{S}}(\bA )$.
\red{Constraints other} than full
symmetry may be relevant in some application areas, such as partial symmetry
as in \textsc{indscal} \cite{CarrC70:psy, TenbSR04:laa}, or
positivity/non-negativity \cite{L1, SmilBG04, TenbKK93:jclassification}.

The definition of symmetric rank is not vacuous because of the following result.

\begin{lemma}
\label{symCanD-prop} Let $\bA \in\mathsf{S}^{k}(\mathbb{C}^{n})$. Then
there exist $\mathbf{y}_{1},\dots,\mathbf{y}_{s}\in\mathbb{C}^{n}$ such that
\[
\bA =\sum\nolimits_{i=1}^{s}\mathbf{y}_{i}^{\otimes k}.
\]
\end{lemma}

\begin{proof}
What we actually have to prove, is that the vector space generated by the
$k$th powers of linear forms $L(\mathbf{x})^{k}$ (for all $L\in\mathbb{C}^{n}$) is not included in a hyperplane of $\mathsf{S}^{k}(\mathbb{C}^{n})$.
This is indeed true, because otherwise there would exist a non-zero element
of $\mathsf{S}^{k}(\mathbb{C}^{n})$ which is orthogonal, \red{under} the bilinear form
$\langle\cdot,\cdot\rangle$, to all $L(\mathbf{x})^{k}$ for $L\in\mathbb{C}^{n}$.
Equivalently, by Lemma \ref{lemDual}, there exists a non-zero polynomial
$q(\mathbf{x})$ of degree $k$ such that $q(L)=0$ for all $L\in\mathbb{C}^{n}$.
But this is impossible, since a non-zero polynomial does not vanish identically on
$\mathbb{C}^{n}$.
\end{proof}

Lemma \ref{symCanD-prop} may be viewed as a particular case of \red{a basic} result
in algebraic geometry, stating that the linear space generated by points of
an algebraic variety that is not included in a {hyperplane, i.e.\ a subspace of codimension $1$,} is the whole space \cite{Harr98, CoxLO98, Shaf77}. For completeness, a proof \red{of our special case is given above}.
Note that it follows from the \red{proof that}
\[
\operatorname*{rank}\nolimits_{\mathsf{S}}(\bA )\leq\binom{n+k-1}{k}
\]
{for all $A \in \mathsf{S}^{k}(\mathbb{C}^{n})$.}

On the other hand, given a symmetric tensor $\bA$, one can compute its
outer product decomposition either in $\mathsf{S}^{k}(\mathbb{C}^{n})$ or in $\mathsf{T}
^{k}(\mathbb{C}^{n})$. Since the outer product decomposition in $\mathsf{S}^{k}
(\mathbb{C}^{n})$ is constrained, it follows \red{that for} all $\bA
\in\mathsf{S}^{k}(\mathbb{C}^{n})$\red{,}
\begin{equation}
\operatorname*{rank}(\bA )\leq\operatorname*{rank}\nolimits_{\mathsf{S}
}(\bA ).\label{rversusrs:eq}%
\end{equation}
We will show that equality holds generically \red{when} $\operatorname*{rank}
_{\mathsf{S}}(\bA )\leq n$ and when $k$ is sufficiently large with
respect to $n$, and always holds {when} $\operatorname*{rank}_{\mathsf{S}%
}(\bA )=1,2$. {While we do not know if the equality holds in general, we suspect that this is the case as we are unaware of any counterexample.}

\subsection{Secant varieties of the Veronese variety}

Let us recall here \red{the correspondence between symmetric outer product decompositions} and secant varieties \red{of the Veronese variety}. By the bijective correspondence between symmetric tensors and homogeneous polynomials established in \eqref{eq:symm-hpoly}, we may discuss this in the context of homogeneous polynomials. The set of homogeneous polynomials {that may be written as a $k$th power of a linear form, $\beta(x)^k=(\beta_{1}x_{1}+\dots+\beta_{n}x_{n})^{k}$ for $\boldsymbol{\beta} =(\beta_{1},\dots,\beta_{n}) \in \mathbb{C}^{n}$,} is a closed algebraic set.
{We may consider this construction as a map from $\mathbb{C}^{n}$ to the space of symmetric tensors given by
\begin{align*}
\nu_{n,k} : \mathbb{C}^{n} & \to \mathbb{C}[x_{1},\dots,x_{n}]_{k} \cong \mathsf{S}^{k}(\mathbb{C}^{n}),\\
 \boldsymbol{\beta} &\mapsto \beta(x)^k.
\end{align*}
The image $\nu_{n,k}(\mathbb{C}^{n})$} is called the \textit{Veronese variety} {and is denoted} $\mathcal{V}_{n,k}$ \cite{Harr98, Z}.
Following this point of view, a {symmetric} tensor is of {symmetric} rank $1$ if it corresponds to a point on the Veronese variety. A {symmetric} tensor is of {symmetric} rank $r$ if it is a linear combination of $r$ {symmetric} tensors of {symmetric} rank $1$ {but not a linear combination of $r - 1$ or fewer such tensors}. In other words, {a symmetric tensor is of symmetric rank not more than $r$ if} it is in the linear space spanned by $r$ points of the Veronese variety.
The closure of the union of all linear spaces spanned by $r$ points of the Veronese variety $\mathcal{V}_{n,k}$ is called\footnote{\red{This seemingly odd choice, i.e.\ $r - 1$ instead of $r$, is standard \cite{Harr98, Z}. The reason being that one wants to be consistent with the usual meaning of a secant, i.e.\ $1$-secant, as a line intersecting \textit{two} points in the variety.}} the \textit{\red{$(r - 1)$th}-secant variety} of $\mathcal{V}_{n,k}$. See \cite{Harr98, Z} for examples and general properties of these algebraic sets. In the {asymmetric} case, {a corresponding notion is obtained by considering the \textit{Segre variety}, i.e.\ the image of the Segre map defined in Section \ref{sec-arrays}.}

\subsection{Why rank can exceed dimension}

We are now in a position to state and prove the following proposition, which is related to a classical result in algebraic geometry stating that $r$ points in $\CC^{n}$ form the solution set
of polynomial equations of degree $\le r$ \cite[pp.\ 6]{Harr98}.
This implies that we can find a  polynomial of degree $\le r-1$ that vanishes
at $r-1$ of the points $L_{i}$ but not at the last one, and hence the independence of polynomials \red{$L_{1}^{r-1},\dots,L_{r}^{r-1}$} follows. Since this proposition is important to our discussion in Section \ref{sec-rank=symrank} (via its corollary below), we give a direct and simple proof \red{below.}

\begin{proposition} \label{propPowerLinear}
Let $L_{1},\dots,L_{r}\in\CC[x_{1},\dots
,x_{n}]_{1}$, i.e.\ linear forms in $n$ variables. If for all
$i\neq j$, $L_{i}$ is not a scalar multiple of $L_{j}$, then for
any $k\geq r-1$, the polynomials \red{$L_{1}^{k},\dots,L_{r}^{k}$} are linearly
independent in $\CC[x_{1},\dots ,x_{n}]$.
\end{proposition}

\begin{proof} Let $k\ge r-1$.
Suppose that for some $\lambda_{1},\dots,\lambda_{r}$,
$\sum\nolimits_{i=1}^{r}\lambda_{i}L_{i}^{k}=0$.
\red{Hence, by} the duality property of Lemma \ref{lemDual}, %
\[
\sum\nolimits_{i=1}^{r}\lambda_{i}\langle F,L_{i}^{k}\rangle=\sum\nolimits_{i=1}^{r}\lambda_{i}F(L_{i}) =0
\]
\red{for all $F\in\CC[x_{1},\dots,x_{n}]_{k}$.}
Let us prove that we can find a homogeneous polynomial $F$ of
degree $k$ \red{that} vanishes at $L_{1},\dots,L_{r-1}$ and
not at $L_{r}$.

Consider a homogeneous polynomial $F$ of degree $k\ge r-1$
\red{that} is a multiple of the product of $r-1$ linear forms $H_{i}$
vanishing at $L_i$ \red{but not at} $L_r$.
We have $F(L_{r})\neq0$ but $F(L_{j})=0$, $1\le j\leq r-1$. As a consequence,
we must have $\lambda_{r}=0$. By a similar argument, we \red{may} show that $\lambda_{i}=0$ for all $i=1,\dots,r$.
\red{It follows} that \red{the} polynomials \red{$L_{1}^{k},\dots,L_{r}^{k}$} are linearly independent.
\end{proof}

Notice that the bound $r-1$ on the degree can be reduced by \red{$d$} if a \red{$d$-dimensional} linear space
containing any \red{$d+1$} of these points does not contain one of the other points \cite[pp.\ 6]{Harr98}red{.}
In this case, we can replace the product of \red{$d+1$} linear forms $H_{i}$ vanishing at \red{$d+1$}
points by \red{just} $1$ linear form vanishing at these \red{$d+1$} points.

\begin{corollary}\label{cor-linIndep}
Let $\bv_1,\dots,\bv_r \in \CC^n$ be $r$ pairwise linearly independent vectors. For any integer {$k\geq r-1$}, the rank-$1$ symmetric tensors
\[
\bv_1^{\otimes k},\dots,\bv_r^{\otimes k} \in \mathsf{S}^{k}(\mathbb{C}^{n})
\]
are linearly independent.
\end{corollary}

This corollary extends results of \cite[Lemma 2.2, pp.\ 2]{Dela06:preprint}
and \cite[Appendix]{KagaLR73}.
Note that vectors
$\mathbf{v}_1,\dots,\mathbf{v}_r$ need not be linearly independent.

\begin{example}
Vectors $\mathbf{v}_{1}=(1,0)$, $\mathbf{v}_{2}=(0,1)$, and $\mathbf{v}%
_{3}=(1,1)$, are pairwise non-collinear but linearly dependent. According to
Corollary \ref{cor-linIndep}, \red{the symmetric tensors $\mathbf{v}_1^{\otimes k}, \mathbf{v}_2^{\otimes k}, \mathbf{v}_3^{\otimes k}$ are linearly independent for any $k \ge 2$}. Evidently, we see that this holds true for
$k=2$ since the matrix below has rank $3$:
\[
\begin{bmatrix}
1 & 0 & 0 & 0\\
0 & 0 & 0 & 1\\
1 & 1 & 1 & 1
\end{bmatrix}.
\]
\end{example}

\subsection{Genericity}

\red{Roughly speaking}, a property is referred to as \textit{typical} if it holds true on a non-zero-volume set {and} \textit{generic} if is true almost everywhere. {Proper} definitions will {follow later} in Section
\ref{sec-main}. It is important to distinguish between typical and generic
properties; for instance, as will be subsequently seen, there can be several
\textit{typical ranks}, but by definition only a single \textit{generic rank}. \red{We will see that} there can be only one typical rank over $\CC$, \red{and} is thus generic.

Through the bijection \eqref{polytens:eq}, the symmetric outer product decomposition
\eqref{cand:eq}
of symmetric tensors can be carried over to quantics, as pointed out in \cite{ComoM96:SP}. The bijection allows one to talk \red{indifferently about} the symmetric outer product decomposition of order-$k$ symmetric tensors \red{and} the decomposition of degree-$k$ quantics into a sum of linear forms raised to the $k$th power.

For a long time, it was believed that there was no explicit expression
for the generic rank. As Reznick pointed out in \cite{Rezn92:mams}, Clebsh proved that even when the numbers of free parameters are the same on both sides of the symmetric outer product decomposition, the generic
rank may not be equal to $\frac{1}{n}\binom{n+k-1}{k}$. For example, in the case
\red{$(k,n)=(4,3)$}, there are $\binom{6}{4} = 15$ degrees of freedom but the generic symmetric rank $R_{\mathsf{S}}(4,3) =6 \neq 5=\frac{1}{3}\binom{6}{4}$. \red{In fact, this} holds true \red{over both $\mathbb{R}$ \cite{Rezn92:mams} and $\mathbb{C}$ \cite{EhreR93}.}
{In Section \ref{sec-values}, we will see} that the generic rank in
$\mathsf{S}^{k}(\mathbb{C}^{n})$ is now known for any order and dimension {due to the ground breaking work of Alexander and Hirschowitz.}

The special case of cubics ($k=3$) is much better known --- a
complete classification is known since 1964 though a constructive algorithm to compute
the symmetric outer product decomposition has only been proposed recently
\cite{KogaM02:issac}. \red{The} simplest case of binary quantics ($n=2$) has also \red{been known} for more than two decades \cite{Wein84:laa, ComoM96:SP, KungR84:ams} --- a result that is used in real world engineering problems \cite{Como04:ieeesp}.


\section{Rank and symmetric rank}\label{sec-rank=symrank}

{Let $\overline{R}_{\mathsf{S}}(k,n)$ be the generic symmetric rank and $R_{\mathsf{S}}(k,n)$ be the maximally attainable symmetric rank} in the space of symmetric tensors $\mathsf{S}^{k}(\mathbb{C}^{n})$. {Similarly, let $\overline{R}(k,n)$ be the generic rank and $R(k,n)$ be the maximally attainable rank in the space of order-$k$ dimension-$n$ cubical tensors $\mathsf{T}^{k}(\mathbb{C}^{n})$. Since $\mathsf{S}^{k}(\mathbb{C}^{n})$ is a subspace of $\mathsf{T}^{k}(\mathbb{C}^{n})$,} generic and maximal ranks (when they exist) are related for every fixed order $k$ and dimension $n$ as follows:
\begin{equation}
\overline{R}(k,n)\geq\overline{R}_{\mathsf{S}}(k,n),\quad\text{and}\quad R(k,n)\geq R_{\mathsf{S}}(k,n). \label{RversusRs:eq}
\end{equation}
It may seem odd that \red{the inequalities in} \eqref{RversusRs:eq} and \eqref{rversusrs:eq} \red{are reversed,}
but there is no contradiction since the spaces are not the same.

It is then legitimate to ask oneself whether the symmetric rank and the rank
are {always equal}. We show that this holds generically when $\operatorname*{rank}%
_{\mathsf{S}}(\bA )\leq n$ (Proposition \ref{propGen1}) or when the
order $k$ is sufficiently large relative to the dimension $n$ (Proposition
\ref{propGen2}). This always holds (not just generically) when
$\operatorname*{rank}_{\mathsf{S}}(\bA )=1,2$ (Proposition
\ref{propGen3}). We will need some preliminary results in proving these assertions.

\begin{lemma}\label{s-linear-indep-lemma}
Let $\mathbf{y}_{1},\dots,\mathbf{y}_{s}\in\mathbb{C}^{n}$ be linearly
independent. Then the symmetric tensor defined by%
\[
\bA :=\sum\nolimits_{i=1}^{s}\mathbf{y}_{i}^{\otimes k}%
\]
has $\operatorname*{rank}_{\mathsf{S}}(\bA )=s$.
\end{lemma}

\begin{proof}
Suppose $\operatorname*{rank}_{\mathsf{S}}(\bA )=r$. Then there exist
$\mathbf{z}_{1},\dots,\mathbf{z}_{r}\in\mathbb{C}^{n}$ such that%
\begin{equation}
\sum\nolimits_{i=1}^{s}\mathbf{y}_{i}^{\otimes k}=\bA =\sum\nolimits_{j=1}^{r}%
\mathbf{z}_{j}^{\otimes k}. \label{eq1}%
\end{equation}
By the linear independence of $\mathbf{y}_{1},\dots,\mathbf{y}_{s}$, there
exist covectors $\varphi_{1},\dots,\varphi_{s}\in(\mathbb{C}^{n})^{\ast}$ that
are dual to $\mathbf{y}_{1},\dots,\mathbf{y}_{s}$, i.e.%
\[
\varphi_{i}(\mathbf{y}_{j})=%
\begin{cases}
1 & \text{if }i=j,\\
0 & \text{if }i\neq j.
\end{cases}
\]
Contracting both sides of \eqref{eq1} in the first $k-1$ modes with
$\varphi_{i}^{\otimes(k-1)}\in\mathsf{S}^{k-1}((\mathbb{C}^{n})^{\ast})$, we
get%
\[
\mathbf{y}_{i}=\sum\nolimits_{j=1}^{r}\alpha_{j}\mathbf{z}_{j},
\]
where $\alpha_{j}=\varphi_{i}(\mathbf{z}_{j})^{k-1}$. In other words,
$\mathbf{y}_{i}\in\operatorname*{span}\{\mathbf{z}_{1},\dots,\mathbf{z}_{r}%
\}$. \red{Since this holds for each $i = 1,\dots, s$, it implies that the $s$ linearly independent vectors $\mathbf{y}_{1},\dots,\mathbf{y}_{s}$ are contained in $\operatorname*{span}\{\mathbf{z}_{1},\dots,\mathbf{z}_{r}\}$. Hence we must have} $r\geq s$. On the other hand, it is clear that $r\leq s$.
\red{Thus we must have equality.}
\end{proof}


\begin{lemma}
\label{genli}Let $s\leq n$. Let $\bA \in\mathsf{S}^{k}(\mathbb{C}^{n})$
with $\operatorname*{rank}_{\mathsf{S}}(\bA )=s$ and%
\[
\bA =\sum\nolimits_{i=1}^{s}\mathbf{y}_{i}^{\otimes k}%
\]
be a symmetric outer product decomposition of $\bA $. Then vectors of the
set $\{\mathbf{y}_{1},\dots,\mathbf{y}_{s}\}$ are generically linearly independent.
\end{lemma}

\begin{proof}
We will write
\[
\mathcal{Y}_{s}:=\{\bA \in\mathsf{S}^{k}(\mathbb{C}^{n})\mid
\operatorname*{rank} \nolimits_{\mathsf{S}}(\bA )\leq s\}\quad
\text{and\quad}\mathcal{Z}_{s} :=\{\bA \in\mathsf{S}^{k}(\mathbb{C}%
^{n})\mid\operatorname*{rank} \nolimits_{\mathsf{S}}(\bA )=s\}.
\]
Define the map from the space of $n\times s$ matrices to order-$k$ symmetric
tensors,
\begin{align*}
f:\mathbb{C}^{n\times s}  &  \rightarrow\mathsf{S}^{k}(\mathbb{C}^{n}),\\
[\mathbf{y}_{1},\dots,\mathbf{y}_{s}]  &  \mapsto\sum\nolimits_{i=1}^{s}
\mathbf{y}_{i}^{\otimes k}.
\end{align*}
It is clear that $f$ takes $\mathbb{C}^{n\times s}$ \textit{onto}
$\mathcal{Y}_{s}$ (i.e.\ $f(\mathbb{C}^{n\times s})=\mathcal{Y}_{s}$). We let
$E_{0}$ and $E_{1}$ be the subsets of rank-deficient and full-rank matrices in
$\mathbb{C}^{n\times s}$ respectively. So we have the disjoint union
\[
{E_{0}\cup E_{1}=\mathbb{C}^{n\times s}, \quad E_{0}\cap E_{1} = \varnothing.}
\]
Recall that the full-rank matrices are generic in $\mathbb{C}^{n\times s}$.
Recall also that $E_{0}$ is an algebraic set in $\mathbb{C}^{n\times s}$
defined by the vanishing of all $s\times s$ principal minors. By the previous
{lemma}, $f(E_{1})\subseteq\mathcal{Z}_{s}$. The set of \red{symmetric tensors,}
\[
\sum\nolimits_{i=1}^{s}\mathbf{y}_{i}^{\otimes k}
\]
in $\mathcal{Z}_{s}$ for which $\{\mathbf{y}_{1},\dots,\mathbf{y}_{s}\}$ is
linearly dependent, i.e.\ $[\mathbf{y}_{1},\dots,\mathbf{y}_{s}]$ is rank
deficient, is simply
\[
\mathcal{Z}_{s}\cap f(E_{0}).
\]
Since $f$ is a polynomial map and $E_{0}$ is a non-trivial algebraic set, we
conclude that $f(E_{1})$ is \red{generic} in $\mathcal{Z}_{s}$.
\end{proof}


\begin{proposition}
\label{propGen1}Let $\bA \in\mathsf{S}^{k}(\mathbb{C}^{n})$. If
$\operatorname*{rank}_{\mathsf{S}}(\bA )\leq n$, then
$\operatorname*{rank}(\bA )=\operatorname*{rank}_{\mathsf{S}}
(\bA )$ generically.
\end{proposition}

\begin{proof}
Let $r=\operatorname*{rank}(\bA )$ and $s=\operatorname*{rank}
_{\mathsf{S}}(\bA )$. So there exist decompositions
\begin{equation}
\sum\nolimits_{j=1}^{r}\mathbf{x}_{j}^{(1)}\otimes\dots\otimes\mathbf{x}_{j}
^{(k)}=\bA =\sum\nolimits_{i=1}^{s}\mathbf{y}_{i}^{\otimes k}.\label{eq2}%
\end{equation}
By Lemma \ref{genli}, we may assume that for a generic $\bA %
\in\mathcal{Z}_{s}$, the vectors $\mathbf{y}_{1},\dots,\mathbf{y}_{s}$ are
linearly independent. As in the proof of Lemma \ref{s-linear-indep-lemma}, we may find a set of covectors $\varphi_{1}%
,\dots,\varphi_{s}\in(\mathbb{C}^{n})^{\ast}$ that are dual to $\mathbf{y}%
_{1},\dots,\mathbf{y}_{s}$, i.e.%
\[
\varphi_{i}(\mathbf{y}_{j})=%
\begin{cases}
1 & \text{if }i=j,\\
0 & \text{if }i\neq j.
\end{cases}
\]
Contracting both sides of \eqref{eq2} in the first $k-1$ modes with
{$\varphi_{i}^{\otimes (k-1)} \in\mathsf{S}^{k-1}((\mathbb{C}%
^{n})^{\ast})$}, we get%
\[
\sum\nolimits_{j=1}^{r}\alpha_{ij}\mathbf{x}_{j}^{(k)}=\mathbf{y}_{i},
\]
where $\alpha_{ij}=\varphi_{i}(\mathbf{x}_{j}^{(1)})\cdots\varphi
_{i}(\mathbf{x}_{j}^{(k-1)})$, $j=1,\dots,r$. Since this holds for each
\red{$i=1,\dots,s$}, it implies that the $s$ linearly independent vectors
$\mathbf{y}_{1},\dots,\mathbf{y}_{s}$ are contained in $\operatorname*{span}%
\{\mathbf{x}_{1}^{(k)},\dots,\mathbf{x}_{r}^{(k)}\}$. Hence we must have
$r\geq s$. On the other hand, it is clear that $r\leq s$. Thus we must have equality.
\end{proof}

\blue{We will see below that we could have $\operatorname*{rank}(\bA)=\operatorname*{rank}_{\mathsf{S}}(\bA )$ even when the constituting vectors $\mathbf{y}_{1},\dots,\mathbf{y}_{s}$ are not linearly independent. The authors would like to thank David Gross for his help in correcting an error in the original proof.}

\begin{proposition}
\label{propGen2}\blue{Let $\mathbf{y}_{1},\dots,\mathbf{y}_{s} \in \mathbb{C}^{n}$ be pairwise linearly independent. If $k$ is sufficiently large, then the symmetric tensor defined by
\[
\bA :=\sum\nolimits_{i=1}^{s}\mathbf{y}_{i}^{\otimes k}%
\]
satisfies $\operatorname*{rank}(\bA )=\operatorname*{rank}_{\mathsf{S}}
(\bA )$ generically.}
\end{proposition}

\begin{proof}
Let $r=\operatorname*{rank}(\bA )$ and $s=\operatorname*{rank}%
_{\mathsf{S}}(\bA )$. So there exist decompositions%
\begin{equation}
\sum\nolimits_{j=1}^{r}\mathbf{x}_{j}^{(1)}\otimes\dots\otimes\mathbf{x}_{j}%
^{(k)}=\bA =\sum\nolimits_{i=1}^{s}\mathbf{y}_{i}^{\otimes k}. \label{eq3}%
\end{equation}
Note that \blue{the \textsc{lhs} may be written $\sum\nolimits_{i=1}^{s}\mathbf{y}_{i}^{\otimes k/2 }\otimes\mathbf{y}_{i}^{\otimes k/2}$, where we have assumed, without loss of generality, that $k$ is even}. By Proposition
\ref{propPowerLinear}, when $k$ is sufficiently large, the
\blue{order-$(k/2)$ tensors $\mathbf{y}_{1}^{\otimes k/2 },\dots,\mathbf{y}%
_{s}^{\otimes k/2 }$} are generically linearly independent. Hence we may find functionals
\blue{$\Phi_{1},\dots,\Phi_{s}\in\mathsf{S}^{k/2}(\mathbb{C}^{n})^{\ast}$} that
are dual to \blue{$\mathbf{y}_{1}^{\otimes k/2},\dots,\mathbf{y}_{s}^{\otimes
 k/2} \in \mathsf{S}^{k/2}(\mathbb{C}^{n})$}, i.e.
\[
\blue{\Phi_{i}(\mathbf{y}_{j}^{\otimes k/2 })} =%
\begin{cases}
1 & \text{if }i=j,\\
0 & \text{if }i\neq j.
\end{cases}
\]
Contracting both sides of \eqref{eq3} in the first \blue{$k/2$} modes with
{$\Phi_{i}$}, we get\blue{
\[
\sum\nolimits_{j=1}^{r}\alpha_{ij}\mathbf{x}_{j}^{(k/2+1)}\otimes \dots \otimes \mathbf{x}_{j}^{(k)}=\mathbf{y}_{i}^{\otimes k/2},
\]}%
where \blue{$\alpha_{ij}=\Phi_{i}(\mathbf{x}_{j}^{(1)}\otimes\dots\otimes
\mathbf{x}_{j}^{(k/2)})$}, $j=1,\dots,r$. Since this holds for each
$i=1,\dots,s$, it implies that the $s$ linearly independent vectors
\blue{$\mathbf{y}_{1}^{\otimes k/2},\dots,\mathbf{y}_{s}^{\otimes k/2}$} are contained in \blue{$\operatorname*{span}%
\{\mathbf{x}_{1}^{(k/2+1)}\otimes \dots \otimes \mathbf{x}_{1}^{(k)},\dots,\mathbf{x}_{r}^{(k/2+1)}\otimes \dots \otimes \mathbf{x}_{r}^{(k)}\}$}. Hence we must have
$r\geq s$. On the other hand, it is clear that $r\leq s$. Thus we must have equality.
\end{proof}

\begin{proposition}
\label{propGen3}Let $\bA \in\mathsf{S}^{k}(\mathbb{C}^{n})$. If
$\operatorname*{rank}_{\mathsf{S}}(\bA )=1$ or $2$, then
$\operatorname*{rank}(\bA )=\operatorname*{rank}_{\mathsf{S}}%
(\bA )$.
\end{proposition}

\begin{proof}
If $\operatorname*{rank}_{\mathsf{S}}(\bA )=1$, then
$\operatorname*{rank}(\bA )=1$ clearly. If $\operatorname*{rank}%
_{\mathsf{S}}(A)=2$, then%
\[
\bA =\mathbf{y}_{1}^{\otimes k}+\mathbf{y}_{2}^{\otimes k}%
\]
for some $\mathbf{y}_{1},\mathbf{y}_{2}\in\mathbb{C}^{n}$. It is clear that
$\mathbf{y}_{1}$ and $\mathbf{y}_{2}$ must be linearly independent or
otherwise $\mathbf{y}_{2}=\alpha\mathbf{y}_{1}$ implies that%
\[
\bA =(\beta\mathbf{y}_{1})^{\otimes k}
\]
for any $\beta=(1+\alpha^{k})^{1/k}$, contradicting $\operatorname*{rank}%
_{\mathsf{S}}(\bA )=2$. It follows from \red{the argument in the proof of} Proposition
\ref{propGen1} with $s=2$ that $\operatorname*{rank}(A)=2$.
\end{proof}

The following result will be useful later.

\begin{proposition}
\label{prop-rank-k} Let $\mathbf{v}_{1}$ and $\mathbf{v}_{2}$ be two linearly
independent vectors in $\mathbb{C}^{n}$. Then for any $k>1$, the following
order-$k$ symmetric tensor:
\begin{multline}
\label{eq-sum-rank-k}\mathbf{v}_{1}\otimes\mathbf{v}_{2}%
\otimes\mathbf{v}_{2}\otimes\dots\otimes\mathbf{v}_{2}+
\mathbf{v}_{2}\otimes\mathbf{v}_{1}\otimes\mathbf{v}%
_{2}\otimes\dots\otimes\mathbf{v}_{2} \\
+ \mathbf{v}_{2}%
\otimes\mathbf{v}_{2}\otimes\mathbf{v}_{1}\otimes\dots
\otimes\mathbf{v}_{2} + \dots+ \mathbf{v}_{2}%
\otimes\mathbf{v}_{2}\otimes\mathbf{v}_{2}\otimes\dots
\otimes\mathbf{v}_{1}%
\end{multline}
is of symmetric rank $k$.
\end{proposition}


\begin{proof}
It is not hard to \red{check that} the symmetric tensor in
\eqref{eq-sum-rank-k} is
associated with the {quantic} {$p(z_1, z_2)=z_{1} z_{2}^{k-1}$},
up to a constant multiplicative factor (where
$z_{1}, z_{2}$ are the first two coordinate variables {in $(z_1,\dots,z_n)$}).

To prove that this {quantic} is of {symmetric} rank $k$, we are going to show that
$p(z_1, z_2)$ can be decomposed into a
sum of powers of linear forms as
\begin{equation}\label{eq:deca}
p(z_1, z_2)=\sum\nolimits_{i=1}^{k} \lambda_{i} (\alpha_{i} z_{1} + \beta_{i}
z_{2})^{k}.
\end{equation}
There are infinitely many possibilities of
choosing coefficients $(\alpha_{i},\beta_{i})$ but {we just} need
to provide one
solution. Take $\alpha_{1}=\dots=\alpha_{r}=1$
and $\beta_{1},\dots, \beta_{k}$ distinct such that\red{
\begin{equation}\label{eq:beta}
\sum\nolimits_{i=1}^{k} \beta_{i}=0.
\end{equation}}%

First we express all quantics in terms of the canonical basis
scaled by the binomial coefficients:
\[
\{z_{1}^{k},k z_{1}^{k-1}z_{2},\dots, k z_{1}z_{2}^{k-1},z_{2}^{k}\}.
\]%
In this basis, the monomial {$k z_{1}%
z_{2}^{k-1}$} can be represented by a $(k+1)$-dimensional vector
containing only
one non-zero entry.
The {quantic} $(z_{i}+\beta_{i}  z_{2})^{k}$  is \red{then} represented by the
vector\red{
\[
[1, \beta_{i}, \beta_{i}^{2}, \dots, \beta_{i}^{k}] \in \mathbb{C}^{k+1}.
\]}%
The existence of coefficients $\lambda_{1}, \dots,\lambda_{k}$
such that we have the decomposition \eqref{eq:deca} is equivalent to the
vanishing of the $(k+1)\times(k+1)$ determinant
\begin{equation}\label{eq:matvdm}
\begin{vmatrix}
0 & 0     & \cdots & 1&  0\\
1 & \beta_{1} &  \cdots & \beta_{1}^{k-1} &\beta_{1}^{k}\\
\vdots & \vdots & & \vdots & \vdots \\
1 & \beta_{k} & \cdots & \beta_{k}^{k-1} & \beta_{k}^{k}\\
\end{vmatrix}.
\end{equation}
An explicit computation shows that this determinant is
$\pm (\sum\nolimits_{i=1}^{k} \beta_{i}) V_{k}(\beta_{1},\dots,\beta_{k})$
where $V_{k}(\beta_{1},\dots,\beta_{k})$ is the {Vandermonde} determinant of
degree $k-1$ of $\beta_{1},\dots,\beta_{k}$.
Thus by \red{\eqref{eq:beta}}, the determinant \red{in} \eqref{eq:matvdm}
vanishes.

This proves that the {symmetric} rank of $z_{1}z_{2}^{k}$ is $\le k$.
Note that the {symmetric} rank cannot be smaller than $k$ because removing any row of
the matrix of \eqref{eq:matvdm} still yields a matrix of rank $k$, if the
$\beta_{i}$ are distinct (see also Proposition \ref{propPowerLinear}).
\end{proof}

This proof is constructive, and gives an algorithm to compute a symmetric outer product decomposition of any binary symmetric tensor of the form \eqref{eq-sum-rank-k}. For
example, the reader can check out that \red{the}
decompositions below \red{may} be obtained this way.

\begin{example}
The {quantics} $48z_{1}^{3}z_{2}$ and $60z_{1}^{4}z_{2}$ {are
associated with the symmetric tensors of maximal rank $A_{31}$ and $A_{41}$ respectively.} {Their} {symmetric outer product decompositions are given by}
\begin{align*}
A _{31} & =8(\mathbf{v}_{1}+\mathbf{v}_{2})^{\otimes
4}-8(\mathbf{v}_{1}-\mathbf{v}_{2})^{\otimes 4} -(\mathbf{v}%
_{1}+2\mathbf{v}_{2})^{\otimes 4}+(\mathbf{v}_{1}-2\mathbf{v}%
_{2})^{\otimes 4},\\
A _{41} & =8(\mathbf{v}_{1}+\mathbf{v}_{2})^{\otimes
5}-8(\mathbf{v}_{1}-\mathbf{v}_{2})^{\otimes 5} -(\mathbf{v}%
_{1}+2\mathbf{v}_{2})^{\otimes 5}+(\mathbf{v}_{1}-2\mathbf{v}%
_{2})^{\otimes 5} +48\mathbf{v}_1^{\otimes5}.
\end{align*}

\end{example}

The maximal symmetric rank achievable by symmetric tensors of order $k$ and
dimension $n=2$ {is $k$, i.e.\ $R_{\mathsf{S}}(k,2)=k$.}
One can say that such \red{symmetric} tensors lie on a tangent line to the \red{Veronese} variety of \red{symmetric rank-$1$} tensors.
In \cite{ComasS01:arxiv}, an algorithm has been proposed to decompose binary forms when their rank is not larger than $k/2$; however, this algorithm would not have found the decompositions above since
the symmetric ranks of $A_{31}$ and $A_{41}$ exceed $4/2$ and $5/2$ respectively.


\section{Generic symmetric rank and typical symmetric ranks}\label{sec-main}

For given order and dimension, define the {following subsets of symmetric tensors $\mathcal{Y}%
_{r}:=\{\bA \in\mathsf{S}^{k}(\mathbb{C}^{n})\mid\operatorname*{rank}%
\nolimits_{\mathsf{S}}(\bA )\leq r\}$ and $\mathcal{Z}_{r}:=\{\bA %
\in\mathsf{S}^{k}(\mathbb{C}^{n})\mid\operatorname*{rank}\nolimits_{\mathsf{S}%
}(\bA )=r\}$.} Also, denote {the corresponding Zariski} {closures} by $\overline{\mathcal{Y}
}_{r}$ and $\overline{\mathcal{Z}}_{r}$ respectively. Recall that the Zariski closure \cite{CoxLO98} of a set $\mathcal{S}$ is simply the smallest variety containing
$\mathcal{S}$. For every $r \in \mathbb{N}$, we clearly have
\[
\mathcal{Y}_{r-1}\cup\mathcal{Z}_{r}=\mathcal{Y}_{r} \quad \text{and} \quad \underbrace{\mathcal{Y}_{1}+\dots+\mathcal{Y}_{1}}_{r\text{ copies}} = \mathcal{Y}_{r}.
\]
The quantities $\overline{R}_{\mathsf{S}}(k,n)$ and $R_{\mathsf{S}}(k,n)$ may now be formally defined by
\[
\overline{R}_{\mathsf{S}}(k,n) := \min\{r\mid\overline{\mathcal{Y}}_{r}=\mathsf{S}^{k}(\mathbb{C}%
^{n})\} \quad \text{\red{and}} \quad
R_{\mathsf{S}}(k,n) := \min\{r\mid\mathcal{Y}_{r}=\mathsf{S}^{k}(\mathbb{C}^{n})\}.
\]
By definition, we have $\overline{R}_{\mathsf{S}}(k,n)\leq R_{\mathsf{S}}(k,n)$. We shall prove in this section that a
generic symmetric rank always exists in $\mathsf{S}^{k}(\mathbb{C}^{n})$\red{, i.e.\ there is an $r$ such that $\overline{\mathcal{Z}}_{r} = \mathsf{S}^{k}(\mathbb{C}^{n})$,} and that it is equal to {$\overline{R}_{\mathsf{S}}(k,n)$, thus justifying our naming $\overline{R}_{\mathsf{S}}(k,n)$ the \textit{generic} symmetric rank in Section \ref{sec-rank=symrank}.}

An integer $r$ is not a \textit{typical rank} if $\mathcal{Z}_{r}$ {has zero}
volume, which means that $\mathcal{Z}_{r}$ is contained in a \red{non-trivial}
closed set. This definition is somewhat unsatisfactory since any mention of `volume' necessarily {involves} a choice of measure, which is really irrelevant here. A better definition {is} as follows.

\begin{definition}
An integer $r$ is a typical rank if $\mathcal{Z}_{r}$ is dense with the
Zariski topology, i.e.\ if $\overline{\mathcal{Z}}_{r}=\mathsf{S}^{k}%
(\mathbb{C}^{n})$. When a typical rank is unique, it may be called
\textit{generic}.
\end{definition}

We used the wording ``typical'' in agreement with previous terminologies \cite{BurgCS97, TenbK99:laa, TenbSR04:laa}.
{Since two dense algebraic sets always intersect over $\CC$}, there can only be one typical rank over $\CC$, and hence {is} generic. In the remainder of this section, we will write $\overline{R}_{\mathsf{S}} = \overline{R}_{\mathsf{S}}(k,n)$ and $R_{\mathsf{S}} = R_{\mathsf{S}}(k,n)$. We can then prove the following.

\begin{proposition}
\label{prop-inclusions} \label{ranks:th} The varieties $\overline{\mathcal{Z}}_{r}$ can be ordered by inclusion as follows. If
\[
r_{1} < r_{2} < \overline{R}_{\mathsf{S}} < r_{3} \leq R_{\mathsf{S}},
\]
then
\[
\overline{\mathcal{Z}}_{r_{1}} \varsubsetneq\overline{\mathcal{Z}}_{r_{2}} \varsubsetneq \overline
{\mathcal{Z}}_{\overline{R}_{\mathsf{S}}} \varsupsetneq\overline{\mathcal{Z}}_{r_{3}}.
\]
\end{proposition}

Before proving this proposition, we first state two preliminary results.
Recall that an algebraic variety is {\textit{irreducible}} if it cannot be decomposed as
the union of proper subvarieties {(cf.\ \cite[pp.\ 51]{Harr98} and \cite[pp.\ 34]{Shaf77})}.
In algebraic geometry, it is known that the secant varieties of any irreducible variety are irreducible.
Nevertheless, we
will give a short proof of the following lemma for the sake of
completeness.

\begin{lemma}
\label{irreducible:lemma} The sets \red{$\overline{\mathcal{Y}}_{r}$, $r\geq1$}, are
irreducible algebraic varieties.
\end{lemma}

\begin{proof}
For $r\ge 1$, the variety $\overline{\mathcal{Y}}_{r}$ is the closure of the image
$\mathcal{Y}_{r}$ of the \red{map
\begin{align*}
\varphi_{r}:\mathbb{C}^{n \times r}  &  \rightarrow\mathsf{S}^{k}(\mathbb{C}^{n}),\\
[\mathbf{u}_{1},\dots,\mathbf{u}_{r}]  &  \mapsto{\sum\nolimits_{i=1}^{r}}%
\mathbf{u}_{i}^{\otimes k}.
\end{align*}}%
Consider now two polynomials $f,g$ such that $fg\equiv0$ on $\overline{\mathcal{Y}%
}_{r}$. As $\overline{\mathcal{Y}}_{r}$ is the Zariski closure of $\mathcal{Y}_{r}%
$, this is equivalent to $fg\equiv0$ on $\mathcal{Y}_{r}$ or
\[
(fg)\circ\varphi_{r}=(f\circ\varphi_{r})(g\circ\varphi_{r})\equiv0.
\]
Thus either $f\equiv0$ or $g\equiv0$ on $\mathcal{Y}_{r}$ or equivalently on
$\overline{\mathcal{Y}}_{r}$, which proves that $\overline{\mathcal{Y}}_{r}$ is an
irreducible variety. For more details on properties of parameterized
varieties, see \cite{CoxLO98}. See also the proof of \cite{Stra83:laa,
BurgCS97} for {third} order tensors.
\end{proof}

\begin{lemma}
\label{minrbar:lemma} We have {$\overline{R}_{\mathsf{S}}=\min\{r \mid\overline{\mathcal{Y}}_{r} =
\overline{\mathcal{Y}}_{r+1}\}$.}
\end{lemma}

\begin{proof}
Suppose that there {exists} {$r<\overline{R}_{\mathsf{S}}$} such that $\overline{\mathcal{Y}}_{r}%
=\overline{\mathcal{Y}}_{r+1}$. Then since $\overline{\mathcal{Y}}_{r} \subseteq \overline{\mathcal{Y}}_{r}+\mathcal{Y}_{1}%
\subseteq \overline{\mathcal{Y}}_{r+1}=\overline{\mathcal{Y}}_{r}$, we have
\[
\overline{\mathcal{Y}}_{r} = \overline{\mathcal{Y}}_{r}+\mathcal{Y}_{1} = \overline{\mathcal{Y}}_{r}+\mathcal{Y}_{1}+\mathcal{Y}_{1}= \dots = \overline{\mathcal{Y}}_{r}+\mathcal{Y}_{1}+\dots+\mathcal{Y}_{1}.
\]
As the sum of $R_{\mathsf{S}}$ copies of $\mathcal{Y}_{1}$ is $\mathsf{S}^{k}(\mathbb{C}^{n})$,
we deduce that $\overline{\mathcal{Y}}_{r}=\mathsf{S}^{k}(\mathbb{C}^{n})$ and \red{thus}
$r\geq\overline{R}_{\mathsf{S}}$, which contradicts our hypothesis. By definition,
$\overline{\mathcal{Y}}_{\overline{R}_{\mathsf{S}}}=\overline{\mathcal{Y}}_{\overline{R}_{\mathsf{S}}+1}=\mathsf{S}%
^{k}(\mathbb{C}^{n})$, which proves the lemma. See also the proof of
\cite{Stra83:laa} \red{for the asymmetric case}.
\end{proof}

We are now in a position to prove {Proposition \ref{prop-inclusions}.}

\textit{Proof of Proposition} \ref{prop-inclusions}. By Lemma
\ref{minrbar:lemma}, we deduce that for {$r<\overline{R}_{\mathsf{S}}$,}
\[
\overline{\mathcal{Y}}_{r}\neq\overline{\mathcal{Y}}_{r+1}.
\]
As $\overline{\mathcal{Y}}_{r}$ is an irreducible variety, we have $\dim
(\overline{\mathcal{Y}}_{r})<\dim(\overline{\mathcal{Y}}_{r+1})$. As {$\mathcal{Y}%
_{r}\cup\mathcal{Z}_{r+1}=\mathcal{Y}_{r+1}$}, we deduce that
\[
\overline{\mathcal{Y}}_{r}\cup\overline{\mathcal{Z}}_{r+1}=\overline{\mathcal{Y}}_{r+1},
\]
which implies by the irreducibility of $\overline{\mathcal{Y}}_{r+1}$, that
$\overline{\mathcal{Z}}_{r+1}=\overline{\mathcal{Y}}_{r+1}$. Consequently, for
{$r_{1}<r_{2}<\overline{R}_{\mathsf{S}}$}, we have
\[
\overline{\mathcal{Z}}_{r_{1}}=\overline{\mathcal{Y}}_{r_{1}}\varsubsetneq
\overline{\mathcal{Z}}_{r_{2}}=\overline{\mathcal{Y}}_{r_{2}}\varsubsetneq
\overline{\mathcal{Z}}_{\overline{R}_{\mathsf{S}}}=\overline{\mathcal{Y}}_{\overline{R}_{\mathsf{S}}}=\mathsf{S}%
^{k}(\mathbb{C}^{n}).
\]
Let us prove now that if {$\overline{R}_{\mathsf{S}}<r_{3}$}, we have $\overline{\mathcal{Z}}_{r_{3}%
}\varsubsetneq\mathsf{S}^{k}(\mathbb{C}^{n})$. Suppose that $\overline{\mathcal{Z}%
}_{r_{3}}=\mathsf{S}^{k}(\mathbb{C}^{n})$, then $\mathcal{Z}_{r_{3}}$ is dense
in $\mathsf{S}^{k}(\mathbb{C}^{n})$ as well as {$\mathcal{Z}_{\overline{R}_{\mathsf{S}}}$} {in} the
Zariski topology. This implies that {$\mathcal{Z}_{r_{3}}\cap\mathcal{Z}%
_{\overline{R}_{\mathsf{S}}}\neq\varnothing$}, which is false because a tensor cannot have two
different ranks. Consequently, we have $\overline{\mathcal{Z}}_{r_{3}}%
\varsubsetneq\mathsf{S}^{k}(\mathbb{C}^{n})$. \endproof

\begin{proposition}
\label{cor-closeness} If {$1 \leq r\leq R_{\mathsf{S}}$}, then $\mathcal{Z}_{r}\neq\overline
{\mathcal{Z}}_{r}$.
\end{proposition}

\begin{proof}
Let $r>1$ and $A \in\mathcal{Z}_{r}$. Then by definition of $\mathcal{Y}_{r}$,
there exists $A_{0} \in\mathcal{Y}_{r-1}$ and $A_{1} \in\mathcal{Y}_{1}$ such
that $A = A_{0} + A_{1}$. As $A_{0} \not \in \mathcal{Y}_{r-2}$ (otherwise
$A\in\mathcal{Y}_{r-1}$) we have $A_{0} \in\mathcal{Z}_{r-1}$. For
$\varepsilon\neq0$, define $A_{\varepsilon} = A_{0} + \varepsilon A_{1}$. We have
that $A_{\varepsilon} \in\mathcal{Z}_{r}$, for all $\varepsilon\neq0$, and
$\lim_{\varepsilon\rightarrow0} A_{\varepsilon} = A_{0}$. This shows that $A_{0}
\in\overline{\mathcal{Z}}_{r} -\mathcal{Z}_{r}$, and consequently that
$\mathcal{Z}_{r}\neq\overline{\mathcal{Z}}_{r}$.
\end{proof}


The above proposition is about the set of symmetric tensors of symmetric rank \textit{exactly} $r$.
But what about those of symmetric rank \textit{at most} $r$? While $\mathcal{Y}_{1}$ is closed as a
determinantal variety, {we will see from Examples \ref{ex-binaryQuantic-closeness} and \ref{ex-ternaryCubic-closeness} as well as Proposition \ref{prop-closeness}} that $\mathcal{Y}_{r}$ is generally not closed for $r>1$.
This is another major difference \red{from} matrices, for which all $\mathcal{Y}%
_{r}$ are closed sets.

\begin{example}
\label{ex-binaryQuantic-closeness} In dimension $n\geq2$, and for any order
$k>2$, $\mathcal{Y}_{2}$ is not closed. In fact, take two independent vectors
$\mathbf{x}_{i}$ and $\mathbf{x}_{j}$ and define the sequence of symmetric tensors{
\begin{equation}
\label{eq-rank2sequence}A _{\varepsilon}(i,j) := \frac{1}{\varepsilon}
\left[  (\mathbf{x}_{i}+\varepsilon\mathbf{x}_{j})^{\otimes k} -
\mathbf{x}_{i}^{\otimes k} \right].
\end{equation}}%
For any \red{$\varepsilon \ne 0$}, $A _{\varepsilon}(i,j)$ is of
symmetric rank $2$, but converges in the limit as $\varepsilon \to 0$ to a {symmetric tensor of symmetric rank $k$}. In fact, the limiting symmetric tensor is easily seen to be a sum of $k$ rank-$1$ \red{tensors,}
\[
\mathbf{x}_{i}\otimes\mathbf{x}_{j}\otimes\dots
\otimes\mathbf{x}_{j} + \mathbf{x}_{j}\otimes\mathbf{x}%
_{i}\otimes\dots\otimes\mathbf{x}_{j} + \dots+
\mathbf{x}_{j}\otimes\mathbf{x}_{j}\otimes\dots
\otimes\mathbf{x}_{i},
\]
which has symmetric rank $k$ by Proposition \ref{prop-rank-k}.
\end{example}

\begin{example}
\label{ex-ternaryCubic-closeness} Let $n=3$ and $k=3$. Then $\mathcal{Y}%
_{5}\subset\overline{\mathcal{Y}}_{3}$, whereas $3<\overline{R}_{\mathsf{S}}$. In fact, take the
symmetric tensor associated with the ternary cubic $p(x,y,z)=x^{2}%
y-xz^{2}$. According to \cite{ComoM96:SP, Rezn94}, this tensor has
rank $5$. On the other hand, it is the limit of the sequence $p_{\varepsilon
}(x,y,z)=x^{2}y-xz^{2}+\varepsilon z^{3}$ as $\varepsilon$ tends to zero.
\red{According to a result in} \cite{ComoM96:SP}, the latter polynomial is associated with a
rank-$3$ tensor since \red{the determinant of} its Hessian is equal to \red{$8x^2(x-3\varepsilon z)$}
and hence contains two distinct linear forms as long as $\varepsilon\neq0$.
\end{example}

It is easy to show that this lack of closeness extends in general to $r>\overline{R}_{\mathsf{S}}$ or
for $r\leq n$, as stated in the two propositions below.

\begin{proposition}
\label{prop-closeness} If {$\overline{R}_{\mathsf{S}}<r$}, then for all $k>2$, $\mathcal{Y}%
_{r}\neq\overline{\mathcal{Y}}_{r}$.
\end{proposition}

\begin{proof}
{If $\overline{R}_{\mathsf{S}}<r$}, then {$\mathcal{Y}_{\overline{R}_{\mathsf{S}}} \varsubsetneq\mathcal{Y}_{r}$}.
{By the} definition \red{of generic} symmetric rank, {$\overline{\mathcal{Y}}_{\overline{R}_{\mathsf{S}}}=\mathsf{S}^{k}(\mathbb{C}^{n})=\overline{\mathcal{Y}}_{r}$}. Hence $\mathcal{Y}%
_{r}\varsubsetneq\overline{\mathcal{Y}}_{r}=\mathsf{S}^{k}(\mathbb{C}^{n})$.
\end{proof}

\begin{proposition}
If $1<r\leq n$, then for any $k>2$, $\mathcal{Y}_{r}\neq\overline{\mathcal{Y}}_{r}$.
\end{proposition}

\begin{proof}
Take $n$ linearly independent vectors $\mathbf{x}_{1},\dots,\mathbf{x}_{n}$. Then
the symmetric tensors $\mathbf{x}_1^{\otimes k},\dots,\mathbf{x}_n^{\otimes k}$
are linearly independent as well, and
$\sum\nolimits_{i=1}^{r} \mathbf{x}_i^{\otimes k}$ is of symmetric rank $r$ for every
$r\leq n$ \red{by Lemma \ref{s-linear-indep-lemma}}. Now for $r>2$ and any $\varepsilon\neq0$, define the symmetric tensor\red{
\[
A _{\varepsilon}= \frac{1}{\varepsilon} \left[(\mathbf{x}_{1}+\varepsilon
\mathbf{x}_{2})^{\otimes k} - \mathbf{x}_{1}^{\otimes
k}\right] + \sum\nolimits_{i=3}^{r} \mathbf{x}^{\otimes k}.
\]}%
$A _{\varepsilon}$ is again of symmetric rank $r$ for every $\varepsilon\neq0$, but tends to a symmetric rank $r+1$
tensor (see also Section \ref{topology-sec}). For $r=2$, the same reasoning
applies with\red{
\[
A _{\varepsilon}= \frac{1}{\varepsilon} \left[(\mathbf{x}_{1}+\varepsilon\mathbf{x}_{2})^{\otimes k}
- \mathbf{x}_{1}^{\otimes k}\right].
\]}%
This shows that $\mathcal{Y}_{r}$ is not closed.
\end{proof}

Based on these two propositions, we conjecture the stronger statement that \red{for order $k > 2$}, the set of symmetric tensors of symmetric rank at most $r$ is never closed, \red{even for $r = n+1,\dots,R_{\mathsf{S}}-1$}.
\begin{conjecture}
\label{conjec-closeness} Assume $k>2$ and $n\geq2$. Then $\mathcal{Y}_{r}%
\neq\overline{\mathcal{Y}}_{r}$ for any $r$ such that $1<r<R_{\mathsf{S}}$.
\end{conjecture}

Up to \red{this point}, our \red{study} has been based on the Zariski topology \cite{Shaf77,
CoxLO98}. However it is useful from \red{a} practical point of view to be able to
apply these results to other topologies, for example\red{,} the Euclidean topology. Since the $\mathcal{Y}_{r}$'s are parameterized and \red{are thus} algebraic constructible sets \cite{Shaf77}, and
since the closure of an algebraic constructible set for the Euclidean topology
and the Zariski topology are the same, the \red{results in this paper} holds true for many
other topologies. We have in particular the following result.

\begin{corollary}
Let $\mu$ be a measure on Borel subsets of $\mathsf{S}^{k}(\mathbb{C}^{n})$
with respect to the Euclidean topology on $\mathsf{S}^{k}(\mathbb{C}^{n})$.
Let {$\overline{R}_{\mathsf{S}}$} be the generic {symmetric} rank in $\mathsf{S}^{k}(\mathbb{C}^{n})$. If
$\mu$ is absolutely continuous with respect to the Lebesgue measure on
$\mathsf{S}^{k}(\mathbb{C}^{n})$, then
\[
\mu(\{\bA \in\mathsf{S}^{k}(\mathbb{C}^{n})\mid\operatorname*{rank}%
\nolimits_{\mathsf{S}}(\bA )\neq\overline{R}_{\mathsf{S}}\})=0.
\]
\end{corollary}
In particular, this corollary tells us that $\mathcal{Z}_{\overline{R}_{\mathsf{S}}}$ is also
dense in $\mathsf{S}^{k}(\mathbb{C}^{n})$ with respect to the Euclidean
topology\red{. It also tells us that} the rank of a tensor whose entries are
drawn randomly according to an absolutely continuous distribution
(e.g.\ Gaussian) is $\overline{R}_{\mathsf{S}}$ with probability $1$. This is useful in signal
processing for instance, where cumulant tensors are estimated from actual
data, and are asymptotically Gaussian distributed \cite{Bril81, Mccu87}.

These statements extend previous results \cite{AtkiL80:laa}, and prove that
there can be only \textit{one} subset $\mathcal{Z}_{r}$ of {non-empty} interior,
and that the latter is dense in $\mathsf{S}^{k}(\mathbb{C}^{n})$; this result{, however,
requires that we work over} an algebraically closed field {such as $\mathbb{C}$}.

The results of this section are indeed not \red{generally} valid over $\mathbb{R}$.
We refer the reader to Section \ref{sec-examples} for further discussions concerning the real field.


\section{Values of the generic symmetric rank\label{sec-values}}

In practice, it would be useful to be able to compute the symmetric rank of any given symmetric
tensor, or at least to know the maximal values of the symmetric rank, given
its order and dimensions. Unfortunately, \red{these questions are far from resolved.}

\red{The corresponding problem for the generic values of the symmetric rank, however, has seen enormous progress due to the work of Alexander and Hirschowitz described in Section \ref{sec:AHT}. In fact, even before their breakthrough,} bounds on the generic symmetric rank have been known for decades
\cite{AtkiL80:laa, Rezn94, Rezn92:mams}:
\[
\left\lceil \frac{1}{n}\binom{n+k-1}{k}\right\rceil \leq\overline{R}_{\mathsf{S}}(k,n)\leq
\binom{n+k-2}{k-1}.
\]
\red{It is known that the} lower bound is often accurate but the
upper bound is not tight \cite{ComoM96:SP}. \red{Furthermore, exact results are
known in the case of binary quantics ($n=2$) and}
ternary cubics ($k=3$) \cite{EhreR93,ComoM96:SP,Rezn92:mams,KogaM02:issac}.

\subsection{\red{Alexander-Hirschowitz Theorem}\label{sec:AHT}} \red{It was not until the work \cite{AlexH95:jag} of Alexander and Hirschowitz in 1995 that the generic symmetric rank problem was completely settled. Nevertheless,} the relevance of their result has remained largely unknown \red{in} the applied and computational mathematics communities.
One reason is that the connection between our problem and the \textit{interpolating polynomials} discussed in \cite{AlexH95:jag} is not \red{at all} well-known in the aforementioned circles. So for the convenience of our readers, we \red{will} state the result of Alexander and Hirschowitz in \red{the context} of the symmetric outer product decomposition below.

\begin{theorem}[Alexander-Hirschowitz]\label{thm-AH}
For $k>2$, the generic {symmetric} rank of {an order-$k$} symmetric tensor
of dimension $n$ {over $\mathbb{C}$} is always equal to the lower bound\red{
\begin{equation}\label{eq:AHT}
\overline{R}_{\mathsf{S}}(k,n) = \left\lceil \frac{1}{n}\binom{n+k-1}{k}\right\rceil
\end{equation}}%
except for the following cases: $(k,n)\in\{(3,5), (4,3), (4,4), (4,5)\}$, {where} it should be increased by $1$.
\end{theorem}

This theorem is extremely complicated to prove, and the interested reader
should refer to the two papers {of Alexander and Hirschowitz} \cite{AlexH95:jag, AlexH92:im}. Simplifications to this proof have also been recently proposed in \cite{Chandler02}. \red{It is worth noting} that these results have been proved in terms of multivariate polynomials and interpolation theory, and not in terms of symmetric tensors.
The exception $(k,n)=(4,3)$ has been known \red{since 1860; in fact, Sylvester referred to it as Clebsh Theorem in his work \cite{Sylv1886:cras}. It is not hard to guess the formula in \eqref{eq:AHT} by a degrees-of-freedom argument. The difficulty of proving Theorem \ref{thm-AH} lies in establishing the fact that the four given exceptions to the expected formula \eqref{eq:AHT} are the only ones.} Table \ref{table-ranks} \red{below lists} a few values of the generic symmetric rank.

\begin{table}[th]
\begin{center}%
\begin{tabular}
[c]{|c||ccccccccc|}\hline
~${}_{k} \;\;{}^{n}$ & ~~2~~ & ~~3~~ & ~~4~~ & ~~5~~ & ~~6~~ & ~~7~~ & ~~8~~ &
~~9~~ & ~~10~~~\\\hline\hline
3 & 2 & 4 & 5 & \textbf{8} & 10 & 12 & 15 & 19 & 22\\
4 & 3 & \textbf{6} & \textbf{10} & \textbf{15} & 21 & 30 & 42 & 55 & 72\\
5 & 3 & 7 & 14 & 26 & 42 & 66 & 99 & 143 & 201\\
6 & 4 & 10 & 21 & 42 & 77 & 132 & 215 & 334 & 501\\\hline
\end{tabular}
\end{center}
\caption{Values of the generic {symmetric rank} $\overline{R}_{\mathsf{S}}(k,n)$ for various
orders $k$ and {dimensions} $n$. Values appearing in bold are the exceptions
outlined by the Alexander-Hirschowitz Theorem.}%
\label{table-ranks}%
\end{table}

\begin{table}[th]
\begin{center}%
\begin{tabular}
[c]{|c||ccccccccc|}\hline
~${}_{k} \;\;{}^{n}$ & ~~2~~ & ~~3~~ & ~~4~~ & ~~5~~ & ~~6~~ & ~~7~~ & ~~8~~ &
~~9~~ & ~~10~~~\\\hline\hline
3 & 0 & 2 & 0 & \textbf{5} & 4 & 0 & 0 & 6 & 0\\
4 & 1 & \textbf{3} & \textbf{5} & \textbf{5} & 0 & 0 & 6 & 0 & 5\\
5 & 0 & 0 & 0 & 4 & 0 & 0 & 0 & 0 & 8\\
6 & 1 & 2 & 0 & 0 & 0 & 0 & 4 & 3 & 5\\\hline
\end{tabular}
\end{center}
\caption{Generic dimension $F(k,n)$ of the fiber of solutions.}%
\label{table-sol}%
\end{table}

\subsection{\red{Uniqueness}}

Besides the exceptions pointed out in Theorem \ref{thm-AH}, the number of
solutions for the {symmetric} outer product decomposition {has to} be finite if the rank $r$
is smaller {than} or equal to {$\frac{1}{n}\binom{n+k-1}{k}$}.
This occurs for instance for all cases of degree $k=5$ in Table
\ref{table-ranks}, except for $n=5$ and $n=10$. {Hence we may deduce the following:}
\begin{corollary}
Suppose $(k,n)\not\in\{(3,5), (4,3), (4,4), (4,5)\}$. Let $A \in\mathsf{S}^{k}(\CC^{n})$ be a generic element and let the symmetric outer product decomposition of $A$ be
\begin{equation}\label{eq:sopd1}
A = \sum\nolimits_{i=1}^{\overline{R}_\mathsf{S}} \mathbf{v}_i^{\otimes k}.
\end{equation}
Then \eqref{eq:sopd1} has a finite number of solutions if and only if
\[
\frac{1}{n} \binom{n+k-1}{k} \in \mathbb{N}.
\]
\end{corollary}


Actually, one \red{may} easily check the generic dimension of the fiber of solutions by computing the number of remaining free parameters \cite{ComoM96:SP}:
\[
F(k,n) = n \overline{R}_{\mathsf{S}}(k,n) - \binom{n+k-1}{k}.
\]
This is summarized in Table \ref{table-sol}. When the dimension of the fiber
is non-zero, there are infinitely many symmetric outer product decompositions.

Our technique is different from the reduction to simplicity proposed by ten
Berge \red{et} al.\ \cite{TenbK99:laa, Tenb04b:jchemo}, but also relies on the
calculation of dimensionality.


\section{Examples\label{sec-examples}}

We will present a few examples to illustrate our discussions in the previous sections.

\subsection{Lack of closeness}\label{topology-sec}

It has been shown \cite{ComoM96:SP, KogaM02:issac} that symmetric tensors of order $3$ and dimension $3$ have a
generic rank {$\overline{R}_{\mathsf{S}}(3,3)=4$} and a maximal rank {$R_{\mathsf{S}}(3,3)=5$}. From the results of
Section \ref{sec-main}, this means that only $\mathcal{Z}_{4}$ is dense in
$\overline{\mathcal{Y}}_{4}=\overline{\mathcal{Y}}_{5}$, and that $\mathcal{Z}_{3}$ and
$\mathcal{Z}_{5}$ are not closed by \red{Proposition} \ref{cor-closeness}. On the
other hand, $\mathcal{Z}_{1}$ is closed.

In order to make this statement even more explicit, let {us} now define a
sequence of {symmetric tensors, each of symmetric rank $2$, that} converges to {a symmetric tensor of symmetric rank $3$}. This will be a simple {demonstration} of the lack of closure of
$\mathcal{Y}_{r}$ for $r>1$ and $k>2$, already stated in Proposition
\ref{prop-closeness}. For this purpose, let {$\mathbf{x},\mathbf{y}$} be two
non-{collinear} vectors. Then the following order-$3$ {symmetric} tensor is of {symmetric} rank $2$ for any
scalar $\varepsilon\neq0$:
\begin{equation}
\label{order3rank2-eq}A _{\varepsilon} = \varepsilon^{2}%
(\mathbf{x}+\varepsilon^{-1}\mathbf{y})^{\otimes3} + \varepsilon
^{2}(\mathbf{x}-\varepsilon^{-1}\mathbf{y})^{\otimes3}
\end{equation}
and it converges, as $\varepsilon \to 0$, to the following symmetric tensor:
\[
A _{0} = 2\left( \mathbf{x}\otimes\mathbf{y}\otimes\mathbf{y}%
+\mathbf{y}\otimes\mathbf{x}\otimes\mathbf{y}+\mathbf{y}\otimes\mathbf{y}%
\otimes\mathbf{x}\right).
\]
This limiting symmetric tensor is of symmetric rank $3$. In fact, one may show \cite{Como02:oxford} that it admits the following symmetric outer product decomposition:
\[
A _{0} = (\mathbf{x}+\mathbf{y})^{\otimes3} - (\mathbf{x}%
-\mathbf{y})^{\otimes3} -2\mathbf{y}^{\otimes3}.
\]
Now let {$\mathbf{x}_{i},\mathbf{y}_{i}$} be linearly independent vectors.

By adding two terms of the form \eqref{order3rank2-eq}, a similar example can
be given in dimension $n=4$, {where we get} a sequence of {symmetric tensors of symmetric rank $4$} converging to a limit of {symmetric} rank $6$.

We will give two more illustrations of Conjecture \ref{conjec-closeness}.

\begin{example}
If the dimension is $n=3$, we can take three linearly independent vectors, say
$\mathbf{x}$, $\mathbf{y}$, and $\mathbf{z}$. Then the sequence of {symmetric} tensors
$A _{\varepsilon}+\mathbf{z}^{\otimes3}$ is of {symmetric} rank $3$ and
converges towards a {symmetric} rank-$4$ tensor.
\end{example}

In dimension $3$, it is somewhat more tricky to build a sequence converging
towards a {symmetric tensor of symmetric rank $5$. Note that $5$ is the maximal rank for $k=3$ and $n=3$.}

\begin{example}
Consider the sequence below as $\varepsilon$ tends to zero:
\begin{equation}
\frac{1}{\varepsilon}\left[  (\mathbf{x}+\varepsilon\mathbf{y}%
)^{\otimes3}-\mathbf{x}^{\otimes3} + (\mathbf{z}%
+\varepsilon\mathbf{x})^{\otimes3}-\mathbf{z}^{\otimes3}
\right].
\end{equation}
It converges to {the following symmetric tensor, which we expressed as a sum of six (asymmetric) rank-$1$ terms,} \[
\mathbf{x}\otimes\mathbf{x}\otimes\mathbf{y}%
+\mathbf{x}\otimes\mathbf{y}\otimes\mathbf{x}%
+\mathbf{y}\otimes\mathbf{x}\otimes\mathbf{x}+
\mathbf{z}\otimes\mathbf{z}\otimes\mathbf{x}%
+\mathbf{z}\otimes\mathbf{x}\otimes\mathbf{z}%
+\mathbf{x}\otimes\mathbf{z}\otimes\mathbf{z}.
\]
{This has symmetric rank $5$ since} it can be associated with {quantic} $x^{2}y+xz^{2}$, which is the
sum of (at least) five cubes.
\end{example}

In terms of algebraic geometry, this example admits a simple geometric interpretation. The limiting tensor is the sum of a point in the tangent space to $\cY_1$ at $\mathbf{x}^{\out3}$ and a point in the tangent space to $\cY_1$ at $\mathbf{z}^{\out3}$.

Note that the same kind of example can be constructed in the asymmetric case:
\[
\mathbf{x}_{1}\otimes\mathbf{x}_{2}\otimes
(\mathbf{x}_{3}-\varepsilon^{-1}\mathbf{y}_{3})+(\mathbf{x}_{1}+\varepsilon
\mathbf{y}_{1})\otimes(\mathbf{x}_{2}+\varepsilon\mathbf{y}_{2})\otimes
\varepsilon^{-1}\mathbf{y}_{3}.
\]
Further discussions of the lack of closeness of $\mathcal{Y}_{r}$ and the ill-posedness of the best rank-$r$ approximation problem \red{in the} asymmetric \red{case} can be found in \cite{dSL}.

\subsection{Symmetric outer product decomposition over the real field}\label{realField-sec}

\red{We now turn our attention to real symmetric tensors. We are interested in
the symmetric outer product decomposition of $A \in \mathsf{S}^{k}(\mathbb{R}^{n})$ over $\mathbb{R}$, i.e.
\begin{equation}\label{eq:rsopd}
A = \sum\nolimits_{i=1}^r \lambda_i \mathbf{v}_i\otimes \mathbf{v}_i\otimes \dots \otimes \mathbf{v}_i
\end{equation}
where $\lambda_i \in \mathbb{R}$ and $\mathbf{v}_i \in \mathbb{R}^n$ for all $i = 1,\dots,r$.
First note that unlike the decomposition over $\mathbb{C}$ in Lemma \ref{symCanD-prop}, we can no longer drop the coefficients $\lambda_1,\dots,\lambda_r$ in \eqref{eq:rsopd} since the $k$th roots of
$\lambda_i$ may not exist in $\mathbb{R}$.}

\red{Since $\mathsf{S}^{k}(\mathbb{R}^{n}) \subset \mathsf{S}^{k}(\mathbb{C}^{n})$, we may regard $A$ as an element of $\mathsf{S}^{k}(\mathbb{C}^{n})$ and seek its symmetric outer product decomposition over $\mathbb{C}$. It is
easy to see that we will generally need more terms in \eqref{eq:rsopd} to decompose $A$ over $\mathbb{R}$ than over $\mathbb{C}$ and so}
\begin{equation}
\operatorname*{rank}\nolimits_{\mathsf{S},\mathbb{C}}(A )\leq\operatorname*{rank}%
\nolimits_{\mathsf{S},\mathbb{R}}(A ). \label{CversusR:eq}%
\end{equation}
This inequality also holds true for the \red{outer product rank} of asymmetric tensors. \red{For $k=2$, i.e.\ matrices, we always have equality in \eqref{CversusR:eq} but we will see in the examples below that strict inequality can occur when $k>2$}.

\begin{example}
\red{Let $A \in \mathsf{S}^{3}(\mathbb{R}^{2}) $ be defined by}
\[
A = \left[
\begin{array}
[c]{rr}%
-1 & 0\\
0 & 1
\end{array}
\right\vert \!\left.
\begin{array}
[c]{rr}%
0 & 1\\
1 & 0
\end{array}
\right].
\]
{It is of symmetric rank $3$ over $\mathbb{R}$}:
\[
A =\frac{1}{2}%
\begin{bmatrix}
1\\
1
\end{bmatrix}
^{\otimes3}+\frac{1}{2}%
\begin{bmatrix}
1\\
-1
\end{bmatrix}
^{\otimes3}-2%
\begin{bmatrix}
1\\
0
\end{bmatrix}
^{\otimes3}%
\]
whereas it is of symmetric rank $2$ over $\mathbb{C}$:
\[
A =\frac{\jmath}{2}%
\begin{bmatrix}
-\jmath\\
1
\end{bmatrix}
^{\otimes3}-\frac{\jmath}{2}%
\begin{bmatrix}
\jmath\\
1
\end{bmatrix}
^{\otimes3}, \quad\text{where }\jmath := \sqrt{-1}.
\]
\red{Hence we see that
$\operatorname*{rank}\nolimits_{\mathsf{S},\mathbb{C}}(A )\ne \operatorname*{rank}\nolimits_{\mathsf{S},\mathbb{R}}(A )$.}
\end{example}

These decompositions \red{may} be obtained \red{using} the algorithm
described in \cite{ComoM96:SP}, for instance. Alternatively, this tensor is
associated with the homogeneous polynomial in two variables $p(x,y)=3xy^{2}-x^{3}$, which can be decomposed over $\mathbb{R}$ into
\[
p(x,y) = \frac{1}{2} (x+y)^{3} + \frac{1}{2} (x-y)^{3} -2x^{3}.
\]

In the case of $2\times2\times2$ symmetric tensors, or equivalently in the
case of binary cubics, the {symmetric} outer product decomposition can always be computed
\cite{ComoM96:SP}. Hence, the {symmetric} rank of any {symmetric} tensor can be calculated, even {over $\mathbb{R}$}. {In this case, it can be shown that the generic symmetric rank over $\mathbb{C}$ is $2$ whereas there are \textit{two typical symmetric ranks} over $\mathbb{R}$, which are $2$ and $3$}

In fact, in the $2\times2\times2$ case, there are two $2\times2$ matrix
slices, that we can call $A_0$ and $A_1$. Since the generic symmetric rank
over $\mathbb{C}$ is $2$, the outer product decomposition is obtained via the eigenvalue
decomposition of the matrix pencil $(A_0, A_1)$, which
generically exists and whose eigenvalues are those of $A_0 A_1^{-1}$. By generating (four) independent real Gaussian entries with zero mean and unit variance, it can be easily checked out with a simple
computer simulation that one gets real eigenvalues in $52$\% of the cases.
This means that the real symmetric rank is $3$ in $48$\% of the remaining
cases. This is the simplest example demonstrating that a generic symmetric rank can be lacking over $\mathbb{R}$. So the concept of typical symmetric rank is essential to studying symmetric tensors over $\mathbb{R}$.

For asymmetric tensors, the same kind of computer simulation would yield (by generating $8$
independent real Gaussian entries) typical ranks of $2$ and $3$, $78$\% and
$22$\% of the time, respectively, leading to the same qualitative conclusions. This procedure is not new
\cite[pp.\ 13]{Tenb04a:jchemo} and has already been proposed in the past to illustrate the existence of several typical ranks for asymmetric tensors \cite{Krus89:mda, TenbK99:laa}. An interesting result obtained by ten Berge \cite{Tenb00:psy} is that $p\times p\times2$ real asymmetric tensors have typical ranks $\{p,p+1\}$.

\red{The problems pertaining to rank and decompositions of real symmetric tensors have not received as much attention as their complex counterparts. However, a moderate amount of work has been done \cite{Krus89:mda, Rezn92:mams, TenbK99:laa, Tenb04b:jchemo, TenbSR04:laa} and we refer the reader to these for further information.}

\subsection{Open questions}

Most of the results that we have presented so far are limited to symmetric tensors
over the complex field. The case of general asymmetric tensors is currently being addressed
with the same kind of approach. As pointed out earlier, decompositions over the
real field are more complicated to handle with algebraic geometric tools. In
addition, while the problem of determining the generic symmetric rank has been resolved thanks to the {Alexander-Hirschowitz Theorem, the maximal symmetric rank is known only for particular values of order and dimensions (e.g.\ dimension $2$); only very rough upper bounds are known for general values. Lastly, the
computation of an explicit symmetric outer product decomposition for a symmetric tensor is computationally expensive, and the conditions (dimension, order) under which this can be
executed within a polynomial time are not yet clearly known. \red{These are problems that we hope will be addressed in future work, either by ourselves or interested readers.}

\subsection*{Acknowledgements}

{The authors would like to thank the anonymous reviewers for
their helpful comments.} This work is a result of collaboration initiated in the 2004
\textit{Workshop on Tensor Decomposition} held at the American
Institute of Mathematics, Palo Alto, CA, and continued in the 2005
\textit{Workshop on Tensor Decomposition and its Applications} held
at the Centre International de Rencontres Math\'{e}matiques (CIRM),
Luminy, France.
The {work} of B.~Mourrain and P.~Comon {has} been partially supported by the contract ANR-06-BLAN-0074 ``\textsc{Decotes}''. The work of G.H.~Golub has been partially supported by the {grant CCF 0430617 from the National Science Foundation}. The work of L.-H.~Lim has been partially supported by the {grant DMS 0101364 from the National Science Foundation}, and by the Gerald J.~Lieberman Fellowship from Stanford University.

\end{document}